\documentclass[11pt]{article}

\usepackage[belowskip=-17pt,aboveskip=5pt]{caption}

\usepackage{amsmath}
\usepackage{amssymb}

\usepackage[english]{babel}

\usepackage[raiselinks=false,plainpages=false,pdfpagelabels,naturalnames=true,colorlinks=true,citecolor=blue,urlcolor=blue,linkcolor=blue,bookmarksopen=true,backref,hyperindex]{hyperref}

\allowdisplaybreaks[1]

\usepackage{latexsym}
\usepackage[emu]{cmll}
\usepackage{prooftree}

\def\Section#1{\section{#1}}
\def\Subsection#1{\subsection{#1}}

\def\Hide#1{\relax}

\def\ANote#1{}
\DeclareSymbolFont{AMSb}{U}{msb}{m}{n}
\DeclareSymbolFontAlphabet{\mathbb}{AMSb}
\DeclareSymbolFont{symbolsC}{U}{txsyc}{m}{n}
\DeclareMathSymbol{\rJoin}{\mathrel}{symbolsC}{89}

\def\DD{\mathbb{D}}

\def\NN{\mathbb{N}}

\def\QQ{\mathbb{Q}}
\def\RR{\mathbb{R}}

\def\ZZ{\mathbb{Z}}

\def\VA{\mathbf{A}}

\def\VC{\mathbf{C}}
\def\VD{\mathbf{D}}
\def\VE{\mathbf{E}}
\def\VF{\mathbf{F}}
\def\VG{\mathbf{G}}
\def\VH{\mathbf{H}}

\def\VK{\mathbf{K}}
\def\VL{\mathbf{L}}
\def\VM{\mathbf{M}}

\def\VP{\mathbf{P}}

\def\VR{\mathbf{R}}
\def\VS{\mathbf{S}}
\def\VT{\mathbf{T}}

\def\VW{\mathbf{W}}

\def\Func#1{{\mathsf{#1}}}

\def\And{\land}
\def\iAnd{\otimes} 

\def\Diff{\mathop{\backslash}}

\def\Dnn{\DD_{{\ge}0}}

\def\Iff{\Leftrightarrow}
\def\I{\Func{I}}
\def\IS{\Func{IS}}

\def\Imp{\Rightarrow}
\def\Inf{\Func{inf}}

\def\IsDef{\mathrel{{:}{=}}}

\def\Ker{\Func{ker}}

\def\Min{\Func{min}}

\def\Not{\lnot}
\def\Or{\lor}

\def\Pmi{\Leftarrow}

\def\ST{\mathrel{|}}

\def\Sup{\Func{sup}}
\def\To{\rightarrow}

\def\Tuple#1{\left\langle #1 \right\rangle}
\def\all#1{\forall #1 {\cdot}\:}
\def\ex#1{\exists #1 {\cdot}\:}

\def\Downset#1{{#1}{\downarrow}}

\newtheorem{Theorem}{Theorem}
\newtheorem{Lemma}[Theorem]{Lemma}
\newtheorem{Corollary}[Theorem]{Corollary}

\newtheorem{Remark}[Theorem]{Remark}

\def\Proof{\par \noindent{\bf Proof: }}
\def\Done{\hfill\rule{0.5em}{0.5em}}
\def\Var{\Func{Var}}

\def\cappedplus#1{\mathop{{+}_{#1}}}

\def\HALVE{[{\sf h}]}
\def\LEASTEL{[{\sf le}]}
\def\Mi{[{\sf m}_i]}
\def\MONE{[{\sf m}_1]}
\def\MTHREE{[{\sf m}_3]}
\def\MTWO{[{\sf m}_2]}
\def\OFOUR{[{\sf o}_4]}
\def\Oj{[{\sf o}_j]}
\def\ONEAX{[{\sf ann}]}
\def\OONE{[{\sf o}_1]}
\def\OTHREE{[{\sf o}_3]}
\def\OTWO{[{\sf o}_2]}
\def\RESIDUAL{[{\sf r}]}

\def\ASM{[\Func{ASM}]}
\def\CE{[{\iAnd}\Func{E}]}
\def\CI{[{\iAnd}\Func{I}]}

\def\Cwc{[{\sf cwc}]}
\def\ATHREE{[{\sf csd}]}

\def\CON{[\Func{CON}]}

\def\CSD{[\Func{CSD}]}
\def\CWC{[\Func{CWC}]}
\def\DN{[\Func{DNE}]}

\def\EFQ{[\Func{EFQ}]}

\def\HL{[\Func{HLB}]}
\def\HU{[\Func{HUB}]}

\def\LE{[{\Lolly}\Func{E}]}
\def\LI{[{\Lolly}\Func{I}]}

\def\Lnot{{{}^{\perp}}}
\def\Lolly{\multimap}

\def\ordSum{\mathop{\stackrel{\frown}{\relax}}}
\def\rImp{\mathop{\rightarrow}}
\def\WK{[\Func{WK}]}

\def\EFQ{[\Func{EFQ}]}

\def\Lang#1{{\cal L}_{\mathbf{#1}}}
\def\Lcl{\Lang{\relax}} 
\def\Lo{\Lang{0}} 
\def\Li{\Lang{1}} 
\def\Lh{\Lang{\frac{1}{2}}} 
\def\Lih{\Lang{1\frac{1}{2}}} 

\def\Logic#1#2{\mbox{{\bf #1}}_{\mbox{\bf #2}}}
\def\ALc{\Logic{AL}{c}}
\def\ALi{\Logic{AL}{i}}
\def\ALu{\Logic{AL}{u}}
\def\BL{\Logic{BL}{\relax}}
\def\CLc{\Logic{CL}{c}}
\def\CLi{\Logic{CL}{i}}
\def\CLu{\Logic{CL}{u}}
\def\IL{\Logic{IL}{\relax}}
\def\ILu{\Logic{IL}{u}}
\def\LLc{\Logic{{\L}L}{c}}
\def\LLi{\Logic{{\L}L}{i}}
\def\LLu{\Logic{{\L}L}{u}}

\title{Hoops, Coops and the Algebraic Semantics of Continuous Logic}

\author{Rob Arthan \& Paulo Oliva}

\begin{document}
\maketitle

\begin{abstract}
B\"{u}chi and Owen studied algebraic structures called hoops. Hoops provide a
natural algebraic semantics for a class of substructural logics that
we think of as intuitionistic analogues of the widely studied {\L}ukasiewicz
logics. Ben Yaacov extended {\L}ukasiewicz logic to get what is called
continuous logic by adding a halving operator. In this paper, we define
the notion of continuous hoop, or coop for short, and show that coops
provide a natural algebraic semantics for continuous logic. We characterise
the simple and subdirectly irreducible coops and investigate the decision
problem for various theories of coops. In passing, we give a new proof
that hoops form a variety by giving an algorithm that converts a
proof in intuitionistic {\L}ukasiwicz logic into a chain of equations.
\end{abstract}

\tableofcontents

\Section{Introduction} \label{sec:introduction}

Around 1930, {\L}ukasiewicz and Tarski \cite{lukasiewicz-tarski30}
instigated the study of logics admitting models in which the
truth values are real numbers drawn from some subset $T$
of the interval $[0, 1]$.
In these models, conjunction is represented by capped addition%
\footnote{We here follow the convention of the literature on continuous logic
in ordering the truth values by increasing logical strength so that $0$
represents truth and $1$ falsehood.}:
$A \And B \IsDef \Inf\{A + B, 1\}$ and
negation is represented by inversion: $\Not A \IsDef 1 - A$. The set $T$
is required to contain $1$ and to be closed under these operations.
One then finds that $T$ is the intersection $G \cap [0, 1]$ where
$G$ is some additive subgroup of $\RR$ with $\ZZ \subseteq G$
and that $T$ is also closed under disjunction and implication defined by
$A \Or B \IsDef \Sup\{A + B - 1, 0\}$ and $A \Imp B \IsDef \Sup\{B - A, 0\}$.
These logics are classical in that $\Not\Not A$ and $A$ are
equivalent. Moreover the law of the excluded middle holds, but
$A \Imp A \And A$ only holds in the special case of Boolean logic
for which $T = \{0, 1\}$, so apart from this special case, the logics
are substructural.

These {\L}ukasiewicz logics have been widely studied, e.g., as instances of fuzzy logics \cite{Hajek98}. More recently ben Yaacov has used
them as a building block in what is called continuous logic
\cite{ben-yaacov-pedersen09}. Continuous logic
unifies work of Henson and others \cite{henson-iovino02} that aims
to overcome shortfalls of classical first-order model theory when
applied to continuous structures such as metric spaces and Banach spaces.
The language of continuous logic extends that of the
usual propositional logic by adding a halving operator, written $A/2$.
In the standard numerical model of continuous logic
the set $T$ of truth values is the interval $[0, 1]$ and
$A/2 \IsDef \frac{1}{2}A$.

Many basic facts about both {\L}ukasiewicz and continuous logics
depend on work of Rose and Rossser \cite{rose-rosser58}
and Chang \cite{chang58a, chang58b} who proved that the following
axiom schemas together with the rule of modus ponens
are complete for the propositional {\L}ukawiecisz logics:
\begin{gather*}
A \Imp (B \Imp A) \tag*{(A1)} \\
(A \Imp B) \Imp (B \Imp C) \Imp (A \Imp C) \tag*{(A2)} \\
((A \Imp B) \Imp B) \Imp ((B \Imp A) \Imp A) \tag*{(A3)} \\
(\Not A \Imp \Not B) \Imp (B \Imp A). \tag*{(A4)}
\end{gather*}
This had been a long-standing conjecture of {\L}ukasiewicz.
Ben Yaacov \cite{ben-yaacov08} added the following axiom schemata
for the halving operator:
\begin{gather*}
(A/2 \Imp A) \Imp A/2 \tag*{(A5)} \\
A/2 \Imp (A/2 \Imp A) \tag*{(A6)}
\end{gather*}
and showed that $A1$--$A6$ together with modus ponens are complete
for the standard numerical model of continuous logic.

The goal of the present paper is to cast some light onto these
axiomatizations by developing propositional {\L}ukasiewicz logic
and continuous logic
as a series of extensions of \emph{intuitionistic} affine logic. A similar
approach for {\L}ukasiewicz logic was developed in \cite{Ciabattoni:1997} where \emph{classical} affine
logic was taken as the starting point. The more restricted setting of intuitionistic affine
logic will allow us to better calibrate the amount of contraction that needs to be added
to affine logic to obtain {\L}ukawiecisz logics. In particular, we obtain an intuitionistic
counter-part of {\L}ukawiecisz logic.

\begin{figure}[t]
\begin{center}
\setlength{\unitlength}{3mm}
\begin{picture}(23,22)(4,7)
\put(4,7){\makebox(3,2){$\ALu$}}
\put(4,17){\makebox(3,2){$\ALi$}}
\put(4,27){\makebox(3,2){$\ALc$}}
\put(14,7){\makebox(3,2){$\LLu$}}
\put(14,17){\makebox(3,2){$\LLi$}}
\put(14,27){\makebox(3,2){$\LLc$}}
\put(20,12){\makebox(3,2){$\ILu$}}
\put(20,22){\makebox(3,2){$\IL$}}
\put(20,32){\makebox(3,2){$\BL$}}
\put(24,7){\makebox(3,2){$\CLu$}}
\put(24,17){\makebox(3,2){$\CLi$}}
\put(24,27){\makebox(3,2){$\CLc$}}
\put(5.5,9.5){\vector(0,1){7.5}}
\put(15.5,9.5){\vector(0,1){7.5}}
\put(25.5,9.5){\vector(0,1){7.5}}
\put(21.5,14.5){\line(0,1){3.25}}
\put(21.5,18.75){\vector(0,1){3.25}}
\put(5.5,19.5){\vector(0,1){7.5}}
\put(15.5,19.5){\vector(0,1){7.5}}
\put(25.5,19.5){\vector(0,1){7.5}}
\put(21.5,24.5){\line(0,1){3.25}}
\put(21.5,28.75){\vector(0,1){3.25}}
\put(7.5,8.25){\vector(1,0){6}}
\put(17.5,8.25){\vector(1,0){6}}
\put(7.5,18.25){\vector(1,0){6}}
\put(17.5,18.25){\vector(1,0){6}}
\put(7.5,28.25){\vector(1,0){6}}
\put(17.5,28.25){\vector(1,0){6}}
\put(17,9){\vector(1,1){3.25}}
\put(17,19){\vector(1,1){3.25}}
\put(17,29){\vector(1,1){3.25}}
\end{picture}
\caption{Relationships between the Logics}
\label{fig:logics}
\end{center}
\end{figure}
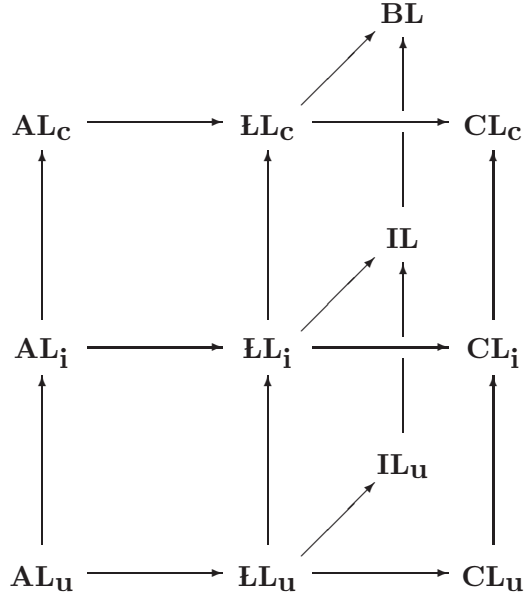

Our work began with the observation that ben Yaacov's continuous logic, which
we call $\CLc$, is an extension of a primitive intuitionistic
substructural logic $\ALi$. We now consider an even more primitive
logic $\ALu$ and develop $\CLc$ as depicted in Figure~\ref{fig:logics},
which also shows how the Brouwer-Heyting intuitionistic propositional logic
$\IL$ and Boolean logic $\BL$ relate to this development.

The structure of the rest of this paper is as follows:
\begin{description}
\item[Section~\ref{sec:the-logics}] gives the definitions of
the logical languages we deal with and of each of the twelve logics
shown in Figure~\ref{fig:logics}.
\item[Section~\ref{sec:algebraic-semantics}] gives sound and complete
algebraic semantics for the logics in terms of certain classes of
pocrims and hoops, algebraic structures that have been quite widely studied in connection
with $\ALi$ and related logics. We introduce the notion of a continuous
hoop or {\em coop} to give the algebraic semantics of continuous logic.
\item[Section~\ref{sec:algebra-of-coops}] considers coops from the
perspective of universal algebra. We characterize the simple and
subdirectly irreducible coops and use the results to begin an
investigation of the decision problem for theories of coops.
\item[Section~\ref{sec:future-work}] outlines further work, particularly
concerning decidability.

\end{description}

\Section{The Logics}\label{sec:the-logics}

We work in a language $\Lcl$ (or $\Lih$ for emphasis) whose
atomic formulas are the propositional constants
$0$ (truth) and $1$ (falsehood) and
propositional variables drawn from
the set $\Var = \{P, Q, \ldots\}$.
If $A$ and $B$ are formulas of $\Lcl$ then so
are $A \iAnd B$ (conjunction),
$A \Lolly B$ (implication) and  $A/2$ (halving).
We define $\Li$ and $\Lh$ to be the sublanguages of $\Lcl$ that disallow halving and $1$ respectively and we define
$\Lo$ to be the intersection of $\Li$ and $\Lh$.
We write $A\Lnot$ as an abbreviation for $A \Lolly 1$.
The judgments of the logics we consider are {\em sequents}
of the form $\Gamma \vdash A$, where the {\em succedent} $A$ is a formula
and the {\em antecedent} $\Gamma$ is a multiset of formulas.

\begin{figure}
\centering
\begin{tabular}{|ccc|}
\hline
\ & \quad \quad & \\
%
%
\begin{prooftree}
\Gamma, A \vdash B
\justifies
\Gamma \vdash A \Lolly B
\using{\LI}
\end{prooftree} & &
%
%
\begin{prooftree}
\Gamma \vdash A \quad \Delta \vdash A \Lolly B
\justifies
\Gamma, \Delta \vdash B
\using{\LE}
\end{prooftree} \\
\ & & \\
%
%
\; \begin{prooftree}
\Gamma \vdash A \quad \Delta \vdash B
\justifies
\Gamma, \Delta \vdash A \iAnd B
\using{\CI}
\end{prooftree} & &
%
%
\begin{prooftree}
\Gamma \vdash A \iAnd B \quad \Delta, A, B \vdash C
\justifies
\Gamma, \Delta \vdash C
\using{\CE}
\end{prooftree} \; \\
\ & & \\
\hline
\end{tabular}
\caption{Inference Rules}
\label{fig:rules}
\end{figure}

The inference rules for all our logics are the
introduction and elimination rules for the two connnectives\footnote{Omitting disjunction from the logic
greatly simplifies the algebraic semantics. While it may be
unsatisfactory from the point of view of intuitionistic philosophy,
disjunction defined using de Morgan's law is adequate for our purposes.} shown in Figure~\ref{fig:rules}.
The various logics we deal with are distinguished only the by axioms we define for them.
We define the axioms such that if  $\Gamma \vdash B$ is an axiom, then so is $\Gamma, A \vdash B$ for any formula $A$.
The way the antecedents of sequents are handled in the inference rules then
implies that we have the following derived rule of weakening:
\[
\begin{prooftree}
\Gamma \vdash B
\justifies
\mathstrut \Gamma, A \vdash B
\using{\WK}
\end{prooftree}
\]
since any instance of this rule in a proof tree may be moved up the proof
tree until it is just beneath an axiom and then the conclusion of the rule
will already be an axiom.

It is easily proved for any logic with the inference rules of
Figure~\ref{fig:rules} that a
form of the deduction theorem holds in the sense that
if one of the following three sequents is provable then so are the other two:
\begin{align*}
A_1, \ldots, A_m &\vdash B, \\
 &\vdash A_1 \Lolly \ldots \Lolly A_m \Lolly B, \\
 &\vdash A_1 \iAnd \ldots \iAnd A_m \Lolly B.
\end{align*}

The axiom schemata for our logics are selected from those shown in Figure~\ref{fig:axioms}.
These are the axiom of assumption $\ASM$, ex-falso-quodlibet $\EFQ$, double negation elimination $\DN$, commutative weak conjunction $\CWC$, commutative strong disjunction $\CSD$, the axiom of contraction $\CON$, and two axioms for the halving operator: one for a lower-bound $\HL$ and one for an upper bound $\HU$.

\begin{figure}
\centering
\begin{tabular}{|c|}\hline
\ \\
%
%
\begin{prooftree}
\justifies
\Gamma, A \vdash A
\using{\ASM}
\end{prooftree}
\ \\[4mm]
%
%
\begin{prooftree}
\justifies
\mathstrut \Gamma, 1 \vdash A
\using{\EFQ}
\end{prooftree}
\  \\[4mm]
%
%
\begin{prooftree}
\justifies
\Gamma, A\Lnot\Lnot \vdash A
\using{\DN}
\end{prooftree}
\ \\[4mm]
%
%
\begin{prooftree}
\justifies
\Gamma, A \iAnd (A \Lolly B) \vdash B \iAnd (B \Lolly A)
\using{\CWC}
\end{prooftree} \;
\  \\[4mm]
%
%
\; \begin{prooftree}
\justifies
\Gamma, (A \Lolly B) \Lolly B \vdash (B \Lolly A) \Lolly A
\using{\CSD}
\end{prooftree}
\ \\[4mm]
%
%
\begin{prooftree}
\justifies
\Gamma, A \vdash A \iAnd A
\using{\CON}
\end{prooftree}
\ \\[4mm]
%
%
\begin{prooftree}
\justifies
\Gamma, A/2 \iAnd A/2 \vdash A
\using{\HL}
\end{prooftree}
\ \\[4mm]
%
%
\begin{prooftree}
\justifies
\Gamma, A/2 \Lolly A \vdash A/2
\using{\HU}
\end{prooftree}
\\\ \\
\hline
\end{tabular}
\caption{Axiom Schemata}
\label{fig:axioms}
\end{figure}

The logics we deal with are discussed in the next few paragraphs
and the axioms for eah logic are summarised in Table~\ref{tab:models}.
In the logics that do not have the axiom schemata $\EFQ$, $1$ plays no special
r\^{o}le and may be omitted from the language and similarly halving
may be omitted from the language in the logics that do not
have the axioms schemata $\HL$ and $\HU$. Alternatively,
the full language $\Lih$ may be used in all cases with
$1$ and $A/2$ effectively acting as variables in formulas that
involve them.

Unbounded\footnote{%
We use the term ``unbounded'' here for logics in which $1$ has no special meaning and need not be an upper bound for the lattice of truth values.
}%
\  intuitionistic affine logic, $\ALu$, has for its axiom schemata $\ASM$ alone. All our other logics include $\ALu$.
Since the contexts $\Gamma$, $\Delta$ are multisets,
an assumption in the rules of Figure~\ref{fig:rules} can be used at most once.
$\ALu$ serves as a prototype for substructural logics with this property.
It corresponds under
the Curry-Howard correspondence to a $\lambda$-calculus with
pairing (i.e., $\lambda$-abstractions of the form $\lambda (x, y)\bullet t$,
$\lambda((x, y), z)\bullet t$, $\lambda(x, (y, z))\bullet t$ etc.) in which
no variable is used twice.
Intuitionistic affine logic, $\ALi$, extends $\ALu$ with the axiom schemata $\EFQ$.

Classical affine logic, $\ALc$, extends $\ALi$ with the axiom schema $\DN$. It can also be viewed as the extension of the so-called multiplicative fragment of Girard's linear logic \cite{Girard(87B)} by allowing weakening and the axiom schema $\EFQ$. We do not define an ``unbounded'' version of $\ALc$ or its extensions, since $\EFQ$ is derivable from $\DN$ in the presence of weakening $\WK$.

Unbounded intuitionistic {\L}ukasiewicz logic, $\LLu$, extends $\ALu$ with
the axiom schema $\CWC$. In $\ALu$, $A \iAnd (A \Lolly B)$
can be viewed as a weak conjunction of $A$ and $B$.
In $\LLu$, we have commutativity of this weak conjunction.
Intuitionistic {\L}ukasiewicz logic, $\LLi$, extends $\LLu$ with the axiom schema $\EFQ$.

Classical {\L}ukasiewicz logic, $\LLc$, extends $\ALi$ with the axiom
schema $\CSD$. In $\ALi$, $(A \Lolly B) \Lolly B$ can be viewed as
a form of disjunction, stronger than that defined by
$(A\Lnot \iAnd B\Lnot)\Lnot$.
In $\ALc$ we have commutativity of this strong disjunction.
This gives us the widely-studied multi-valued logic of {\L}ukasiewicz.

\begin{table}
\centering
\begin{tabular}{|c|l|l|} \hline
{\bf Logic} & {\bf Axioms} & {\bf Models} \\
\hline \hline
$\ALu$ & $\ASM$ & pocrims \\[1mm]
\hline
$\ALi$ & as $\ALu + \EFQ$ & bounded pocrims \\[1mm]
\hline
$\ALc$ & as $\ALi + \DN$& involutive pocrims \\[1mm]
\hline
$\LLu$ & as $\ALu + \CWC$ & hoops \\[1mm]
\hline
$\LLi$ & as $\LLu + \EFQ$ & bounded hoops \\[1mm]
\hline
$\LLc$ & as $\ALi + \CSD$ & bounded involutive hoops \\[1mm]
\hline
$\ILu$  & as $\ALu + \CON$ & idempotent pocrims \\[1mm]
\hline
$\IL$  & as $\ALi + \CON$ & bounded idempotent pocrims \\[1mm]
\hline
$\BL$  & as $\IL + \DN$ & involutive idempotent pocrims \\[1mm]
\hline
$\CLu$ & as $\LLu + \HL + \HU$ & coops \\[1mm]
\hline
$\CLi$ & as $\LLi + \HL + \HU$ & bounded coops \\[1mm]
\hline
$\CLc$ & as $\LLc + \HL + \HU$ & involutive coops \\[1mm]
\hline
\end{tabular}
\caption{The logics and their models}
\label{tab:models}
\end{table}

Unbounded intuitionistic propositional logic, $\ILu$, extends $\ALu$ with the axiom
schema $\CON$,
which is equivalent to a contraction rule allowing $\Gamma, A \vdash B$
to be derived from $\Gamma, A, A \vdash B$.
Intuitionistic propositional logic, $\IL,$ extends $\ILu$ with the axioms schemata $\EFQ$.
$\IL$ is the conjunction-implication fragment of the well-known Brouwer-Heyting intuitionistic propositional logic.

Boolean logic, $\BL$, extends $\IL$ with the axiom schema $\DN$.
This is the familiar two-valued logic of truth tables.

Unbounded intuitionistic continuous logic, $\CLu$, allows the halving operator
and extends $\LLu$ with the
axiom schemas $\HL$ and $\HU$, which effectively give lower and upper bounds
on the logical strength of $A/2$. The two axioms can also be read as saying that
$A/2$ is equivalent to $A/2 \Lolly A$.
Intuitionistic continuous logic, $\CLi$, extends $\CLu$ with the axiom schemata $\EFQ$.
$\CLi$ is an intuitionistic version of the continuous
logic of ben Yaacov \cite{ben-yaacov-pedersen09}.

Classical Continuous logic, $\CLc$ extends $\CLi$ with the axiom
schema $\DN$. This gives ben Yaacov's continuous logic. The motivating
model takes truth values to be real numbers between $0$ and $1$ with
conjunction defined as capped addition.

Our initial goal was to understand the relationships amongst
$\ALi$, $\LLc$ and $\CLc$. The other logics came into focus
when we tried to decompose the somewhat intractable axiom
$\CSD$ into a combination of $\DN$ and an intuitionistic component.
We will see that the twelve logics are related as shown in
Figure~\ref{fig:logics}. In the figure,
an arrow from $T_1$ to $T_2$ means that $T_2$
extends $T_1$, i.e., the set of provable sequents of $T_2$ contains that of
$T_1$. In each square, the north-east logic is the least extension of
the south-west logic that contains the other two. For human beings, at
least, the proof of this fact is quite tricky for the $\ALi$-$\LLc$ square,
see~\cite[chapters 2 and 3]{Hajek98}.

The routes in Figure~\ref{fig:logics} from $\ALu$ to $\IL$ and $\BL$ have been
quite extensively studied \cite{blok-ferreirim00,raftery07}.  We are not aware
of any work on $\CLu$ and $\CLi$, but these are clearly natural objects of
study in connection with ben Yaacov's continuous logic.  It should be noted
that $\ILu$ and $\CLu$ are incompatible: given $\CON$, $A/2$ and $A/2 \iAnd
A/2$ are equivalent, so that from $\HL$ and $\HU$ one finds that $A/2 \Lolly A$
and $A/2$ are both provable; which proves $A$, for arbitrary formulas $A$.

\Section{Algebraic Semantics}\label{sec:algebraic-semantics}

\ANote{We need to explain that this is a fairly potted account and give
references for more detail of the kind of thing we are doing. It would be
nice to have a reference (or give our own account, if necessary) contrasting
our natural deduction/Gentzen style approach with the Hilbert-style approach
in much of the literature.}

We give an algebraic semantics to our logics using pocrims:
partially ordered, commutative, residuated, integral monoids.
A {\em pocrim}\footnote{Strictly speaking, this is a {\em dual pocrim},
since we order it by increasing logical strength and write it additively,
whereas in much of the literature the opposite order and multiplicative
notation is used (so halves would be square roots).
We follow the ordering convention of the continuous logic literature.},
is a structure for the signature $(0, +, \rImp)$ of type
$(0, 2, 2)$ satisfying the following laws:
\begin{gather*}
(x + y) + z = x + (y + z)                               \tag*{$\MONE$} \\
x + y = y + x                                           \tag*{$\MTWO$} \\
x + 0 = x                                             \tag*{$\MTHREE$} \\
x \ge x                                                 \tag*{$\OONE$} \\
x \ge y \And y \ge z \Imp x \ge z                       \tag*{$\OTWO$} \\
x \ge y \And y \ge x \Imp x = y                       \tag*{$\OTHREE$} \\
x \ge y \Imp x + z \ge y + z                           \tag*{$\OFOUR$} \\
x \ge 0                                              \tag*{$\LEASTEL$} \\
x + y \ge z \Iff x \ge y \rImp z                    \tag*{$\RESIDUAL$}
\end{gather*}

\noindent where $x \ge y$ is an abbreviation for $x \rImp y = 0$.

When working in a pocrim, we adopt the convention that $\rImp$ associates
to the right and has lower precedence than $+$. So, for example,
the brackets in $(a + b) \rImp (c \rImp (d + f))$ are all redundant, while
those in $((a \rImp b) \rImp c) + d$ are all required.

Let $VM = (M, 0, +, \rImp)$ be a pocrim.
The laws~$\Mi$, $\Oj$ and $\LEASTEL$ say that
$(M, 0, +; {\ge})$ is a partially ordered commutative monoid
with the identity $0$ as least element.
In particular, $+$ is monotonic in both its arguments.
Law~$\RESIDUAL$, the {\em residuation property},
says that for any $x$ and $z$ the
set $\{y \ST x + y \ge z\}$ is non-empty and has $x \rImp z$
as least element.
As is easily verified, $\rImp$ is antimonotonic in its first
argument and monotonic in its second argument.

Let $\alpha : \Var \To M$ be an interpretation of logical variables
as elements of $M$
and extend $\alpha$ to a function $v_{\alpha} : \Lo \To M$
by interpreting $0$, $\iAnd$ and $\Lolly$ as $0$, $+$ and
$\rImp$ respectively.
We say that $\alpha$ {\em satisfies} the sequent
$C_1, \ldots, C_n \vdash A$ iff $v_{\alpha}(C_1) + \ldots + v_{\alpha}(C_n) \ge v_{\alpha}(A)$.
We say that $\Gamma \vdash A$ is {\em valid} in $\VM$
iff it is satisfied by every assignment $\alpha : \Var \To M$, in which
case we say $\VM$ is a model of $\Gamma \vdash A$.
If $\cal C$ is a class of pocrims, we say $\Gamma \vdash A$ is valid
in $\cal C$ if it is valid in every member of $\cal C$.
We say that a logic $L$ whose language is $\Lo$ 
is sound for a class of pocrims $\cal C$
if every sequent over $\Lo$ that is provable in $L$ is valid in $\cal C$.
We say that $L$ is complete for $\cal C$ if the converse holds.
We then have:

\begin{Theorem}\label{thm:alu-sound-complete}
$\ALu$ is sound and complete for the class of all pocrims.
\end{Theorem}
\Proof
This is standard.
Soundness is a routine exercise. For the completeness, one defines
an equivalence relation $\simeq$ on formulas such that $A \simeq B$ holds
iff both $A \vdash B$ and $B \vdash A$ are provable in the logic.
Writing $[A]$ for the equivalence class of a formula $A$, one then
shows that the set of equivalence classes $T$ is the carrier set
of a pocrim $\VT = (T; 0, +, \rImp)$, where $0 = [0]$ and the
operators $+$ and $\rImp$ are defined so that
$[A] + [B] = [A \iAnd B]$ $[A] \rImp [B] = [A \Lolly B]$.
In $\VT$, the {\em term model} of the logic, $C_1, \ldots, C_n \vdash A$ is
valid, i.e., $[C_1] + \ldots [C_n] \rImp [A] = 0$ holds,
iff $C_1, \ldots, C_n \vdash A$ is provable. Completeness follows, since
a sequent that is valid in all pocrims must be valid in the pocrim $\VT$
and hence must be provable.
\Done

\vspace{4mm}

The above theorem says that a sequent
is provable in $\ALu$ iff it has every pocrim as a model.  In the sequel we
will often use the theorem to derive laws that hold in all pocrims. For
example, it is easy to find a proof in $\ALu$ of the sequent $P, Q \Lolly P
\vdash Q \Lolly (P \iAnd P)$, from which we may conclude that the law $x + (y
\rImp x) \ge y \rImp x + x$ holds in any pocrim.

A {\em hoop} is a pocrim
that is {\em naturally ordered}, i.e., whenever $x \ge y$, there is $z$
such that $x = y + z$. It is a nice exercise in  the use of the residuation
property to show that a pocrim is a hoop iff it satisfies the identity
\[
\begin{array}{l@{\quad\quad}r}
x + (x \rImp y) = y + (y \rImp x)                          \tag*{$\Cwc$}
\end{array}
\]
From this it follows that the logic $\LLu$ is sound and complete for the
class of all hoops.
See~\cite{blok-ferreirim00} for more information on hoops.

We say a pocrim is {\em idempotent}
if it is idempotent as a monoid, i.e., it satisfies $x + x = x$.
Note that this condition implies condition $\Cwc$, since it implies
$x + y \ge x + (x \rImp y) = x + (x + (x \rImp y)) \ge x + y$, whence,
$x + (x \rImp y) = x + y = y + x = y + (y \rImp x)$.
Using this, we find that $\ILu$ is sound and complete for the
class of all idempotent pocrims.

To complete our treatment of the bottom layer in Figure~\ref{fig:logics},
we need to prove a lemma about hoops that will help us with the
algebraic semantics of the halving operator. The hoop axiom $\Cwc$ is
surprisingly powerful but often requires considerable ingenuity to apply.
We have been greatly assisted in our work by using the late Bill McCune's
{\tt Prover9} and {\tt Mace4} programs to prove algebraic facts and to find
counter-examples.
Readers who enjoy a challenge may like to look for their own proof of
the following lemma before reading ours, which is a fairly direct
translation of that found after a few minutes by {\tt Prover9}.
\begin{Theorem}[{\tt Prover9}] \label{thm:hoop-halving-bounds}
The following hold in any hoop:
\begin{align*}
\mbox{\em(i)}\quad&
  \mbox{if $a = a \rImp b$ and $c = c \rImp b$, then $a = c$} \\
\mbox{\em(ii)}\quad&
  \mbox{if $a \ge a \rImp b$ and $c = c \rImp b$, then $a \ge c$} \\
\mbox{\em(iii)}\quad&
  \mbox{if $a \le a \rImp b$ and $c = c \rImp b$, then $a \le c$}.
\end{align*}
\end{Theorem}
\Proof
\noindent
{\em(i):} this is immediate from parts~{\em(ii)} and~{\em(iii)}.

\noindent
{\em(ii):} by the hypothesis on $a$, 
$a + a \ge b$, and so, using the hypothesis on $c$
and the fact that $x + (y \rImp x) \ge y \rImp x + x$
discussed in the remarks following the proof of
Theorem~\ref{thm:alu-sound-complete}, we find:
\[
a + (c \rImp a) \ge c \rImp a + a \ge c \rImp b = c
\]
\noindent so that $c \rImp a \ge a \rImp c$. Using the fact
that if $x \ge y$, then $x + (y \rImp z) \ge z$, we have:
$
(c \rImp a) + ((a \rImp c) \rImp c) \ge c.
$
As the hypothesis on $c$ implies $c + c \ge b$, this gives:
\[
c + (c \rImp a) + ((a \rImp c) \rImp c) \ge b.
\]
\noindent Using $\Cwc$ twice and the fact that $x \rImp y \rImp x = 0$, we find:
\begin{align*}
b &\le c + (c \rImp a) + ((a \rImp c) \rImp c) \\
   &= a + (a \rImp c) + ((a \rImp c) \rImp c) \\
   &= a + c + (c \rImp a \rImp c) \\
   &= a + c.
\end{align*}
\noindent I.e., $a + c \ge b$, so that $a \ge c \rImp b = c$ as required.

\noindent
{\em(iii):} by the hypothesis on $c$ and using
the fact that $(x \rImp y) + x \ge y$ twice, we have: 
\[
c + (a \rImp c) + a = (c \rImp b) + (a \rImp c) + a \ge b.
\]
\noindent Using the hypothesis on $a$, we have:
\[
c + (a \rImp c) \ge a \rImp b \ge a
\]
So $c \ge (a \rImp c) \rImp a$, implying:
\[
c \rImp (a \rImp c) \rImp a = 0.
\]
\noindent Using $\Cwc$ and the facts that $x \rImp y \le (z \rImp x) \rImp y$ and $(x \rImp y) + x \ge y$, we have:
\begin{align*}
c &= c + (c \rImp (a \rImp c) \rImp a) \\
   &= ((a \rImp c) \rImp a) + (((a \rImp c) \rImp a) \rImp c) \\
   &\ge ((a \rImp c) \rImp a) + (a \rImp c) \\
   &\ge a. \tag*{\Done}
\end{align*}

We define a {\em coop} to be a structure for the signature $(0, +, \rImp, /2)$
of type $(0, 2, 2, 1)$ whose $(0, +, \rImp)$-reduct is a hoop and
such that for every $x$ we have:
\[
\begin{array}{l@{\quad\quad}r}
x/2 = x/2 \rImp x                                             \tag*{$\HALVE$}
\end{array}
\]
From $\HALVE$, one clearly has $x \ge x/2 \rImp x = x/2$, i.e., $x \rImp x/2 =
0$ and so using also $\Cwc$ one finds $x/2 + x/2 = x/2 + (x/2 \rImp x) = x + (x
\rImp x/2) = x$ justifying the choice of notation. (Our convention is that $/2$
binds
tighter than the infix operators, so the brackets are needed in $(x + y)/2$
but not in $x \rImp (x/2)$).

The following very useful theorem shows that
the halving operator is uniquely defined by the condition $x/2 = x/2 \rImp x$
(so we could have defined a coop to be a hoop that
satisfies the axiom $\all{x}\ex{y}y = y \rImp x$ and taken the
halving operator to be defined on such a hoop by equation~$\HALVE$).

\begin{Theorem} \label{thm:coop-halving-bounds}
Let $a$ and $b$ be elements of a coop. Then the following hold:
\begin{align*}
\mbox{\em(i)}\quad&
   a = b/2 \Iff a = a \rImp b \\
\mbox{\em(ii)}\quad&
   a \ge b/2 \Iff a \ge a \rImp b \\
\mbox{\em(iii)}\quad&
   a \le b/2 \Iff a \le a \rImp b
\end{align*}
\end{Theorem}
\Proof
$\Imp$: let $R \in \{=, \ge, \le\}$, then using the definition of a coop
and the fact that $\rImp$ is antimonotonic in its left argument, we have:
\[
a \mathrel{R} b/2 = b/2 \rImp b \mathrel{R} a \rImp b
\]
$\Pmi$: immediate from
Theorem~\ref{thm:hoop-halving-bounds} and the definition of a coop. \Done

\begin{Corollary}\label{cor:halving-miscellany}
Let $a$ and $b$ be elements of a coop. Then the following hold:
\begin{align*}
\mbox{\em(i)}\quad&
   a = b \Iff a/2 = b/2 \\
\mbox{\em(ii)}\quad&
   a \ge b \Iff a/2 \ge b/2 \\
\mbox{\em(iii)}\quad&
   a/2 = a \Iff a = 0 \\
\mbox{\em(iv)}\quad&
   a/2 + b/2 \ge (a + b)/2 \\
\mbox{\em(v)}\quad&
   a/2 \rImp b/2 = (a \rImp b)/2
\end{align*}
\end{Corollary}
\Proof
{\em(i):} immediate from {\em(ii)}.
\\
{\em(ii)$\Imp$:} if $a \ge b$, then as $a = a/2 + a/2$, we have $a/2 \ge a/2
\rImp b$ and then, by the theorem, $a/2 \ge b/2$.
\\
{\em(ii)$\Pmi$:} if $a/2 \ge b/2$, then $a = a/2 + a/2 \ge b/2 + b/2 = b$.
\\{\em{(iii)}}: if $a/2 = a$, then, by the theorem,
$a/2 = a/2 \rImp a = a \rImp a = 0$.
\\{\em{(iv)}}: We have
\begin{align*}
a/2 + b/2 \rImp a/2 + b/2 \rImp a + b &= a/2 + b/2 + a/2 + b/2 \rImp a + b \\
   &= a + b \rImp a + b\\
   &= 0
\end{align*}
I.e., $a/2 + b/2 \ge a/2 + b/2 \rImp a + b$,
so, by the theorem, $a/2 + b/2 \ge (a +b)/2$
\\{\em{(v)}}:
I claim that
$(a \rImp b)/2 \ge a/2 \rImp b/2$ and
$(a \rImp b)/2 \le a/2 \rImp b/2$, from which the result follows.
For the first part of the claim, we have $(a \rImp b)/2 \ge a/2 \rImp b/2$ iff
$a/2 + (a \rImp b)/2 \ge b/2$ and, by the theorem, this holds iff
$a/2 + (a \rImp b)/2 \ge a/2 + (a \rImp b)/2 \rImp b$, i.e.
iff
$b \le (a/2 + (a \rImp b)/2) + (a/2 + (a \rImp b)/2) = a + (a \rImp b)$
which is true.
For the second part of the claim, we have:
\begin{align*}
(a/2 \rImp b/2) \rImp a \rImp b
   &= (a/2 \rImp b/2) + a \rImp b \\
   &= a/2 + a/2 + (a/2 \rImp b/2) \rImp b \\
   &= a/2 + b/2 + (b/2 \rImp a/2) \rImp b \tag*{$\Cwc$} \\
   &= a/2 + (b/2 \rImp a/2) \rImp b/2 \rImp b \\
   &= a/2 + (b/2 \rImp a/2) \rImp b/2 \tag*{\HALVE}\\
   &\le a/2 \rImp b/2.
\end{align*}
where the inequality follows from the fact that $\rImp$ is antimonotonic in
its first argument.
So, by the theorem, $(a \rImp b)/2 \le a/2 \rImp b/2$ as required.
\Done
\begin{Corollary}\label{cor:no-finite-coops}
There are no non-trivial finite coops.
\end{Corollary}
\Proof
If $a$ is a non-zero element of a coop, parts {\em(ii)} and {\em(iii)} of the corollary
imply that the sequence $a, a/2, (a/2)/2, \ldots$ is strictly decreasing.
Hence a finite coop has no non-zero elements.
\Done

\vspace{4mm}

Given an interpretation, $\alpha : \Var \To M$ with values in
a coop, we extend the function $v_{\alpha} : \Lo \To M$ to
$\Lh$ in such a way that $v_{\alpha}(A/2) = (v_{\alpha}(A))/2$ and extend the
notions of satisfaction, etc. accordingly.
The proof of Theorem~\ref{thm:alu-sound-complete} is easily extended to show that
the logic $\CLu$ is sound and complete for the class of coops (using
Theorem~\ref{thm:coop-halving-bounds} to show that the halving operation
on the term model is well-defined).

We now have a sound and complete algebraic semantics for each of the logics in
the bottom layer of Figure~\ref{fig:logics}. Moving to the middle layer, let us
say that a pocrim, hoop or coop is {\em bounded} if it has a (necessarily
unique) {\em annihilator}, i.e., an element $1$ such that for every $x$ we have:
\[
\begin{array}{l@{\quad\quad}r}
x + 1 = 1 & \ONEAX
\end{array}
\]
Assume the pocrim $\VM$ is bounded. Then $0 \le x \le x + 1 = 1$ for any $x$ and
$(M; \le)$ is indeed a bounded ordered set.
Given an interpretation, $\alpha : \Var \To M$ with values in
a bounded pocrim, we extend the function $v_{\alpha} : \Lo \To M$
to $\Li$ so that $v_{\alpha}(1) = 1$ and extend the notions of satisfaction
etc. accordingly. Yet 
again the proof of Theorem~\ref{thm:alu-sound-complete} is easily extended to show that
the logic $\ALi$ is sound and complete for the class of bounded pocrims.

We then find that the logics $\LLi$, $\CLi$ and $\IL$ are sound and complete
for bounded hoops, bounded coops and idempotent bounded hoops respectively.
Idempotent bounded hoops are also known as Brouwerian algebras and are known
to be the conjunction-implication reducts of Heyting algebras (see \cite{Koehler81} and the works cited therein).

Finally, for the top layer of Figure~\ref{fig:logics}, we say a
pocrim is {\em involutive} if it is bounded and satisfies $\Not\Not x = x$, where we write
$\Not x$ as an abbreviation for $x \rImp 1$,
Idempotent involutive hoops are easily seen to be the conjunction-implication reducts of Boolean algebras.
We find that $\ALc$, $\LLc$, $\CLc$ and $\BL$ are sound and complete
for involutive pocrims, involutive hoops, involutive coops and
idempotent involutive hoops respectively.
This completes the proof of the following theorem:
\begin{Theorem}\label{thm:all-sound-complete}
The logics $\ALi, \ldots, \CLc$, $\ILu$, $\IL$ and $\BL$
of Figure~\ref{fig:logics}
are sound and complete for the corresponding classes of pocrims,
hoops and coops listed in Table~\ref{tab:models}.
\Done
\end{Theorem}

A {\em Wajsberg hoop} is a hoop satisfying the identity
\[
\begin{array}{l@{\quad\quad}r}
(x \rImp y) \rImp y = (y \rImp x) \rImp x & \ATHREE
\end{array}
\]
It can be shown that Wajsberg hoops are the same as bounded
involutive hoops.

The classes of pocrims associated with the logics in the left-hand column
in Figure~\ref{fig:logics} are very general: any
partial order can be embedded in an involutive pocrim.
To see this, let $X$ be any partially ordered set.
Take a disjoint copy $X\Lnot$ of $X$ (say $X\Lnot = X \times \{1\}$) and write
$x\Lnot$ for the image in $X\Lnot$ of $x \in X$.  Choose objects  $0$, $1$, $r$
and $s$ distinct from each other and from the elements of $X \cup
X\Lnot$ and order the disjoint union $P_X = \{0, r\} \cup X \cup X\Lnot
\cup \{s, 1\}$ so that,
{\em(i),} $0 < r < X < X\Lnot < s < 1$,
{\em(ii),} the subset $X$ has the given ordering and,
{\em(iii),} $X\Lnot$ has the opposite ordering.
Extend the mapping $(\cdot)\Lnot:X \To X\Lnot$ to all of
$P_X$ so that $0\Lnot = 1$, $r\Lnot = s$ and $a\Lnot\Lnot = a$ for all
$a$.  Then $\Lnot$ is an order-reversing mapping of $P_X$ onto itself and there
is a unique commutative binary operation $+$ on $P_X$ with the following
properties:
$$
\begin{array}{rcl@{\quad}l}
a + 0 &=& a, &\mbox{for every $a$;}\\
a + b &=& s, &\mbox{for every $a, b \ge r$ such that $a \not\ge b\Lnot$;}\\
a + b &=& 1, &\mbox{for every $a, b \ge r$ such that $a \ge b\Lnot$.}
\end{array}
$$
Now let $\VP_X = (P_X, 0, +, \rImp)$ where $\rImp$ is defined
using de Morgan's law: $a \rImp b = (a + b\Lnot)\Lnot$.
Then one finds that $a \rImp b = 0$ iff $a \ge b$ in $P_X$ with respect to the
order defined above and the laws for an involutive pocrim other than
associativity of $+$ are then easily verified for $\VP_X$.
For the associativity of $+$, first note
that if $0 \in \{a, b, c\}$, $(a + b) + c = a + (b
+ c)$ is trivial.  If $a, b, c \ge r$ then $a + b, b + c
\ge r\Lnot$ and we have:
\[
1 \ge (a + b) + c \ge r\Lnot + r =  1 = r + r\Lnot \le a + (b + c) \le 1.
\]
\noindent so that $a + (b + c) = 1 = (a + b) + c$. Thus $\VP_X$ is indeed
an involutive pocrim.

It is known that the class of involutive pocrims is not a variety i.e., it
cannot be characterised by equational laws.  Since involutive pocrims are
characterised over bounded pocrims and over pocrims by equational laws, it
follows that the class of pocrims and the class of bounded pocrims are also not
varieties.  See~\cite{raftery07} and the works cited therein for these results
and their history and for further information about pocrims in general and
involutive pocrims in particular.

Bosbach~\cite{bosbach69a} gave a direct proof of an equational axiomatization of
the class of hoops. Using Theorem~\ref{thm:all-sound-complete}, we can give an
alternative proof that shows how a proof of a sequent $\vdash A$ may be
translated into an equational proof that $\alpha = 0$, where $\alpha$ is a
translation into the language of pocrims of the formula $A$.

\begin{Theorem}\label{thm:luk-provability-equational}
A structure $\VH = (H; 0, +, \rImp)$ is a hoop
iff $(H; 0, +)$ is a commutative monoid and the following equations hold in $H$:
\begin{enumerate}
\item\label{luk-imp-self} $x \rImp x = 0$
\item\label{luk-imp-zero} $x \rImp 0 = 0$
\item\label{luk-zero-imp} $0 \rImp x = x$
\item\label{luk-conj-imp} $x + y \rImp z = x \rImp y \rImp z$
\item\label{luk-cwc} $x + (x \rImp y) = y + (y \rImp x)$
\end{enumerate}
\end{Theorem}
\Proof
It follows easily from the definitions (or from
Theorem~\ref{thm:all-sound-complete}) that the equations hold in any hoop.
For the converse, Theorem~\ref{thm:all-sound-complete}
implies that it is sufficient to show that if there is
proof of $\vdash A$ in $\LLu$
then $[A]$ (the element of the term model of $\LLu$
represented by $A$) can be reduced to 0 using the commutative
monoid laws and equations~\ref{luk-imp-self} to \ref{luk-cwc}.
More generally, if $B_1, \ldots, B_m$ and $A$ are formulas, with
$\gamma = [B_1] + \ldots + [B_m]$ and
$a = [A]$, we will show how to translate
a proof of $B_1, \ldots, B_m \vdash A$ into a sequence
of equations $\gamma \rImp a = a_1 = \ldots = a_n = 0$
where each equation $a_i = a_{i+1}$ is obtained
by applying one of the equations~\ref{luk-imp-self} to \ref{luk-cwc}
to a subterm of $a_i$ or $a_{i+1}$ and then simplifying or
rearranging as necessary using the commutative monoid laws.
We have base cases for the axioms $\ASM$ and $\CWC$ of Figure~\ref{fig:axioms}
and inductive steps for the rules of Figure~\ref{fig:rules}.

\noindent
$\ASM$: we want $\gamma + a \rImp a = 0$ for arbitrary $\gamma$ and $a$:
\begin{align*}
\gamma + a \rImp a &=
  & \tag*{(eq. \ref{luk-conj-imp})} \\
 \gamma \rImp a \rImp a &=
  & \tag*{(eq. \ref{luk-imp-self})} \\
 \gamma \rImp 0 &= 0
  & \tag*{(eq. \ref{luk-imp-zero})}
\end{align*}

\noindent
$\CWC$: we want $\gamma + a + (a \rImp b) \rImp b + (b \rImp a) = 0$
for arbitrary $\gamma$, $a$ and $b$:
\begin{align*}
\gamma + a + (a \rImp b) \rImp b + (b \rImp a) &=
   &\tag*{(eq. \ref{luk-cwc})} \\
 \gamma + b + (b \rImp a) \rImp b + (b \rImp a) &=
   &\tag*{(eq. \ref{luk-conj-imp})} \\
 \gamma \rImp b + (b \rImp a) \rImp b + (b \rImp a) &=
   &\tag*{(eq. \ref{luk-imp-self})} \\
 \gamma \rImp 0 &= 0  
   &\tag*{(eq. \ref{luk-imp-zero})}
\end{align*}

\noindent
$\LI$: we are given $\gamma + a \rImp b = 0$ and we want
$\gamma \rImp a \rImp b = 0$:
\begin{align*}
\gamma \rImp a \rImp b &=
     &\tag*{(eq. \ref{luk-conj-imp})}\\
   \gamma + a \rImp b &= 0
    &\tag*{(hyp.)}
\end{align*}

\noindent
$\LE$: we are given $\gamma \rImp a = 0$ and $\delta \rImp a \rImp b =
0$ and we want $\gamma + \delta \rImp b = 0$:
\begin{align*}
\gamma + \delta \rImp b &=
     &\tag*{(hyp.)}\\
   \gamma + (\gamma \rImp a) + \delta \rImp b &=
     &\tag*{(eq. \ref{luk-cwc})}\\
   a + (a \rImp \gamma) + \delta \rImp b &=
     &\relax{}\\
   (a \rImp \gamma) + \delta + a \rImp b &=
     &\tag*{(hyp.)}\\
   (a \rImp \gamma) + \delta + a + (\delta \rImp a \rImp b) \rImp b &=
     &\tag*{(eq. \ref{luk-conj-imp})}\\
   (a \rImp \gamma) + \delta + a + (\delta + a \rImp b) \rImp b &=
     &\tag*{(eq. \ref{luk-cwc})}\\
   (a \rImp \gamma) + b + (b \rImp \delta + a) \rImp b &=
     &\relax{}\\
   (a \rImp \gamma) + (b \rImp \delta + a) + b \rImp b &=
     &\tag*{(eq. \ref{luk-conj-imp})}\\
   (a \rImp \gamma) + (b \rImp \delta + a) \rImp b \rImp b &=
     &\tag*{(eq. \ref{luk-imp-self})}\\
   (a \rImp \gamma) + (b \rImp \delta + a) \rImp 0 &= 0
     &\tag*{(eq. \ref{luk-imp-zero})}
\end{align*}

\noindent
$\CI$: we are given $\gamma \rImp a = 0$ and $\delta \rImp b = 0$
and we want $\gamma + \delta \rImp a + b = 0$.
\begin{align*}
\gamma + \delta \rImp a + b &=
     &\tag*{(hyp.)}\\
   \gamma + (\gamma \rImp a) + \delta + (\delta \rImp b) \rImp a + b &=
     &\tag*{(eq. \ref{luk-cwc})}\\
   a + (a \rImp \gamma) + b + (b \rImp \delta) \rImp a + b &=
     &\relax{}\\
   (a \rImp \gamma) + (b \rImp \delta) + a + b \rImp a + b &=
     &\tag*{(eq. \ref{luk-conj-imp})}\\
   (a \rImp \gamma) + (b \rImp \delta) \rImp a + b \rImp a + b &=
     &\tag*{(eq. \ref{luk-imp-self})}\\
   (a \rImp \gamma) + (b \rImp \delta) \rImp 0 &= 0
     &\tag*{(eq. \ref{luk-imp-zero})}
\end{align*}

\noindent
$\CE$: we are given $\gamma \rImp a + b = 0$ and $\delta + a + b \rImp c =
0$ and we want $\gamma + \delta \rImp c = 0$:
\begin{align*}
\gamma + \delta \rImp c &=
     &\tag*{(hyp.)}\\
   \gamma + (\gamma \rImp a + b) + \delta \rImp c &=
     &\tag*{(eq. \ref{luk-cwc})}\\
   a + b + (a + b \rImp \gamma) + \delta \rImp c &=
    &\relax{} \\
   (a + b \rImp \gamma) + \delta + a + b \rImp c &=
     &\tag*{(eq. \ref{luk-conj-imp})}\\
   (a + b \rImp \gamma) \rImp \delta + a + b \rImp c &=
     &\tag*{(hyp.)}\\
   (a + b \rImp \gamma) \rImp 0 &= 0
     &\tag*{(eq. \ref{luk-imp-zero})} 
\end{align*}
This completes the induction. \Done

\vspace{4mm}

The axiomatization in the statement of Theorem~\
\ref{thm:luk-provability-equational} is natural and convenient but by no means
minimal. See~\cite{bosbach69a} for more concise axiomatizations.

\Section{Algebra of Coops}\label{sec:algebra-of-coops}

Blok and Ferreirim \cite{blok-ferreirim00} 
have studied hoops from the perspective of universal algebra.
Here we undertake an analogous study of coops.
Our goal is to obtain decision problems for useful theories of coops.
This will require
various facts about hoops, most of which may be found in
\cite{blok-ferreirim00}, but in the dual (multiplicative) notation.
We begin by looking at some special classes of coops, for which certain
facts that hold for involutive hoops can be obtained rather efficiently
by dint of the halving operator.

\Subsection{Some Special Classes of Coops}\label{sec:special-classes}

We say a hoop is {\em cancellative} if its underlying monoid is a cancellation
monoid ($x + y = x + z$ implies $y = z$).  Let us say a hoop is {\em
semi-cancellative} if $x + y = x + z$ and $y \not= z$ implies $x + y$ is an
annihilator (i.e., the hoop is bounded with $x + y = 1$).  Thus 
a hoop that is semi-cancellative and not bounded is cancellative.
In a linearly ordered hoop, the semi-cancellative property is easily seen to be
equivalent to the condition that $x + y = x$ implies that either $y = 0$ or $x$
is an annihilator.

Semi-cancellative coops enjoy the property that halving is almost
a homomorphism, or, indeed, a real homomorphism if the coop is cancellative:
\begin{Lemma}\label{lma:semi-cancellative-half-plus}
Let $\VC$ be a semi-cancellative coop, then, for any $x, y \in C$, either
$ (x + y)/2 = x/2 + y/2 $ or $x + y = 1$.
\end{Lemma}
\Proof
Since $x/2 + x/2 = x$ and $y/2 + y/2 = y$, we have
$x + y = x/2 + y/2 + x/2 + y/2$.
On the other hand, since $x + y \ge x/2 + y/2$, $\Cwc$ implies that
$x + y = x/2 + y/2 + (x/2 + y/2 \rImp x + y)$.
By the semi-cancellative property, either $x + y = 1$ or
$x/2 + y/2 = x/2 + y/2 \rImp x + y$.
In the latter case, Theorem~\ref{thm:coop-halving-bounds} (i) tells us that
$x/2 + y/2 = (x + y)/2$.
\Done

\vspace{4mm}

We now prove a very useful theorem that will let us transfer some
important results about bounded coops to unbounded coops. This corresponds
to Chang's construction of the enveloping group of an MV-algebra but
the proof involves much less tricky algebra. 

Before stating the theorem, we introduce some notation and terminology
that will be used throughout the sequel.
Let $\VG = (G; 0, +, \ge)$
be a 2-divisible linearly ordered commutative group.
Writing $G_{{\ge}0}$ for the set of non-negative elements of $G$, we then
have a coop
$\VG_{{\ge}0} = (G_{{\ge}0}; 0, +, \rImp, /2)$ where $x \rImp y \IsDef \Sup\{0, y - x\}$
and $x/2$ is that element of $G$ such that $x/2 + x/2 = x$ (this is unique
because $\VG$ is linearly ordered and hence torsion-free).

If $\VL$ is any coop and $a$ is any non-zero element of $\VL$,
we have a bounded coop $\VL_a = (\{x\in L \ST x \le a\}, \cappedplus{a}, \rImp, /2)$
where $x \cappedplus{a} y \IsDef \Inf\{a, x + y\}$.
We say $\VL_a$ is $\VL$ {\em capped at $a$}.
We will just write $x + y$
for $x \cappedplus{a} y$ in contexts where it is clear that we are working in
$\VL_a$. If $\VL = \VG_{{\ge}0}$ for some
2-divisible linearly ordered commutative group $\VG$, we write
$\VG_{[0, a]}$ for $\VL_a$.
Note that $\VG_{[0, a]}$ is an involutive coop: with $\Not x = a - x$,
we clearly have $\Not\Not x = x$.

As an example, take $\VG$ to be the additive group $\DD$ of dyadic rationals
$\DD = (\{\frac{i}{2^n} \ST i \in \ZZ, n \in \NN\}; 0, +, \ge)$. We then have
an unbounded coop $\Dnn$ and from $\Dnn$, we obtain the bounded coops $\DD_{[0,
a]} = ([0, a] \cap \DD, 0, \cappedplus{a}, \rImp, /2)$ for $a$ any positive
dyadic rational.  Note that the isomorphism type of $\DD_{[0, a]}$ depends on
$a$: $\DD_{[0, 1]}$ contains no $x$ such that $3x$ is the annihilator but
$\DD_{[0, 3]}$ does.

\begin{Theorem}\label{thm:c-hat}
Let $\VC$ be a semi-cancellative bounded coop. Then there exist a
cancellative unbounded
coop $\hat{\VC}$, an element $\hat{1} \in \hat{C}$ and an isomorphism $\alpha :
\VC \To \hat{\VC}_{\hat{1}}$.  Every element of $\hat{\VC}$ has the form
$2^m\alpha(a)$ for some $a \in C$ and $m \in \NN$.  If $\VC$ is linearly ordered
then so is $\hat{\VC}$.
\end{Theorem}
\Proof
Let $\VD = \VC^{\NN}$ be the product of countably many copies of $\VC$.
Thus elements of $\VD$ are sequences
$x = \Tuple{x_0, x_1, \ldots}$ of elements of $C$ and the coop operations
are defined pointwise: $(x + y)_i = x_i + y_i$,
$(x \rImp y)_i = x_i \rImp y_i$ and $(x/2)_i = x_i/2$.
For this proof, let us say $x \in D$ is {\em regular} if
$x_{i+1} = x_i/2$ for all but finitely many $i$.
Using Corollary~\ref{cor:halving-miscellany} and Lemma \ref{lma:semi-cancellative-half-plus} as appropriate, it is easy to see that if $x$ and $y$ are regular then so
are $x \rImp y$, $x + y$ and $x/2$. Thus the regular elements comprise
a subcoop $\VR$ of $\VD$.
Define a relation $\sim$ on $R$ by $x \sim y$
iff $x_i = y_i$ for all but finitely many $i$.
It is a routine exercise to verify that $\sim$ is a congruence.
Let $\hat{\VC}$ be $\VR/{\sim}$ and, for $a \in C$, let
$\alpha(a)$ be given by $(\alpha(a))_i = \frac{1}{2^i}a$ and let
$\hat{1} = \alpha(1)$.
By Corollary~\ref{cor:halving-miscellany}, $\alpha(a \rImp b) = \alpha(a) \rImp \alpha(b)$
for any $a, b \in C$, and, by Lemma \ref{lma:semi-cancellative-half-plus},
if $a + b < 1$, $\alpha(a + b) = \alpha(a) + \alpha(b)$.
It is easy to verify that $\alpha$ is an injection and that
$\alpha(C) = \{a \in \hat{C} \ST \hat{1} \ge a\}$, from which
it follows that $\alpha$ is an isomorphism between $\VC$ and $\hat{\VC}_{\hat{1}}$.
If $x \in R$, there is $a \in C$ and $m \in \NN$ such that
for all $i \in \NN$, $x_{m+i} = \frac{1}{2^i}a$ and then
$[x] = 2^m\alpha(a)$. 
Hence, for any $s \in \hat{C}$, $\frac{1}{2^i}s \in \alpha(C)$
for all but finitely many $i$ and
from this it follows that, $\hat{\VC}$ is semi-cancellative and
hence cancellative and that, if $\VC$ is linearly ordered,
then so is $\hat{\VC}$.
\Done

\begin{Theorem}\label{thm:c-bar}
Let $\VC$ be a linearly ordered cancellative unbounded coop. Then there exist a 2-divisible linearly ordered
group $\overline{\VC}$ and an isomorphism $\beta : \VC \To \overline{\VC}_{{\ge}0}$.
\end{Theorem}
\Proof
Define $\overline{\VC}$ to be the group of differences of $\VC$
and let $\beta : \VC \To \overline{\VC}$ be the natural homomorphism.
Every element of $\overline{\VC}$ has the form $\beta(a) - \beta(b)$ for
$a, b \in C$.
$\beta(a) - \beta(b) = \beta(c) - \beta(d)$ iff there are $x, y \in C$, such
that $a + x = c + y$ and $b + x = d + y$.
We have $(\beta(a/2) - \beta(b/2)) + (\beta(a/2) - \beta(b/2)) =
\beta(a) - \beta(b)$, so $\overline{\VC}$ is 2-divisible.
As $\VC$ is linearly ordered, given $a, b \in C$, either {\em(i)} $a \ge b$,
in which case, $\beta(a) - \beta(b) = \beta(b \rImp a)$,
since $a + 0 = (b \rImp a) + b$ and $b + 0 = 0 + b$, or {\em(ii)}
$b \ge a $, in which case, $\beta(a) - \beta(b) = -\beta(a \rImp b)$,
since $a + 0 = 0 + a$ and $b + 0 = (a \rImp b) + a$.
Thus for any $s \in \overline{C}$, either $s \in \beta(C)$ or $s \in -\beta(C)$.
Moreover if $s \in \beta(C) \cap -\beta(C)$, we have $s = \beta(a) = -\beta(b)$
whence for some $x, y \in C$ we have $a + x = y$ and $x = b + y$, whence
$a + b + y = y$ implying $a = b = 0$, thus $\beta(C) \cap -\beta(C) = \{0\}$.
Since $\beta(C) + \beta(C) = \beta(C)$,
it follows that $\beta(C)$ is the non-negative cone of a linear order on
$\overline{\VC}$ and that $\beta$ is an isomorphism of $\VC$ with 
$\overline{\VC}_{{\ge}0}$.
\Done

\begin{Theorem}\label{thm:lin-canc-coops-decidable}
The first order theories of the following classes of coops are decidable:
{\em(i)} linearly ordered cancellative coops
{\em(ii)} linearly ordered bounded semi-cancellative coops,
{\em(iii)} linearly ordered semi-cancellative coops.
\end{Theorem}
\Proof
Using Theorems~\ref{thm:c-hat} and~\ref{thm:c-bar}, one can find
primitive recursive reductions
of the theory of linearly ordered bounded semi-cancellative coops
to that of linearly ordered cancellative coops and of the latter
theory to the theory of linearly ordered 2-divisible groups.
The theory of linearly ordered
groups is decidable by a well-known result of Gurevich, and hence so
is the theory of 2-divisible linearly ordered groups (since the latter is
a finitely axiomatisable extension of the former).
Hence, {\em(i)} and {\em(ii)} hold.
As for {\em(iii)}, a general linearly ordered
semi-cancellative coop is either cancellative or bounded,
so the theory in {\em(iii)} is the intersection of the theories in {\em(i)}
and {\em(ii)}.
\Done

\Subsection{Homomorphisms and Ideals}

Let $\VH$ be a hoop.
An {\em ideal} $I$ of $\VH$, is a downwards-closed submonoid:
\begin{gather*}
0 \in I \subseteq H \\
I + I \subseteq I \\
\Downset{I} \subseteq I
\end{gather*}
\noindent where, for any $X, Y \subseteq H$,
$X + Y = \{ x + y \ST x \in X, y \in Y\}$
and
$\Downset{X} = \{y \in H \ST \ex{x \in X} x \ge y\}$.
For example, if $X \subseteq H$, the {\em ideal generated by $X$}, $\I(X)$, is
the set comprising all $y \in H$, such that for some $x_1, \ldots, x_n \in X$,
$y \le x_1 + \ldots + x_n$. $\I(X)$ is easily seen to be an ideal and is
clearly the smallest ideal containing $X$. As a special case, the ideal $\I(x)$
generated by $x \in H$, comprises all elements $y$ such that $y \le nx$ for
some $n \in \NN$.  We say an ideal $I$ is {\em proper} if ${0} \not= I \not= H$.

If $I$ is an ideal, then $I$ is actually the carrier set of a subhoop, since,
we have $I \rImp I \subseteq H \rImp I \subseteq I$ (since $I$ is
downwards-closed and $x \rImp y \le y$ for any $x$ and $y$).
If $\VK$ is also a hoop and $f : H \To K$ is a homomorphism of hoops,
we define the {\em kernel} of $f$, $\Ker(f)$, as follows:
\[
\Ker(f) \IsDef \{x : H \ST f(x) = 0\}.
\]
\noindent $\Ker(f)$ is clearly a submonoid of $\VH$. Moreover, if $y \in
\Ker(f)$ and $x \le y$, then, by definition, $f(y) = 0$ and $y \rImp x = 0$,
and then $f(x) = f(y) \rImp f(x) = f(y \rImp x) = f(0) = 0$, so $x \in
\Ker(f)$.  Thus $\Ker(f)$ is an ideal of $\VH$.  Conversely, if $I$ is an ideal
of $\VH$, define a relation $\theta \subseteq H \times H$, by $x
\mathrel{\theta} y \Iff x \rImp y \in I \And y \rImp x \in I$.  It is then
routine to verify that $\theta$ is a hoop congruence on $\VH$ and
that, with $p_{\theta} : \VH \To \VH/\theta$, the natural projection onto
the quotient hoop, we have $\Ker(p_{\theta}) = I$.
It follows that the lattice of congruences on $\VH$ is isomorphic to its
lattice of ideals. In particular, a hoop is simple (i.e., it admits no
non-trivial congruences) iff it has no proper ideals (so that
$\I(x) = H$ for every non-zero $x \in H$).

\begin{Theorem}\label{thm:coop-homomorphisms}
If $\VC$ and $\VD$ are coops then a mapping $f : C \To D$ is a homomorphism
of coops iff it is a homomorphism of the underlying hoops of $\VC$ and $\VD$.
\end{Theorem}
\Proof
Necessity is trivial. For sufficiency,
assume $f: C \To D$ is a homomorphism of hoops. By definition,
$f(x \rImp y) =
f(x) \rImp f(y)$ for any $x, y \in C$.  So for any $x \in C$, we have:
\[
f(x/2) \rImp f(x) = f(x/2 \rImp x) = f(x/2)
\]
\noindent whence by Theorem~\ref{thm:coop-halving-bounds} we must have $f(x/2)
= f(x)/2$. It follows that $f$ is a homomorphism of coops.
\Done

\vspace{4mm}

Thus we need no new notion for the kernels of coop homomorphisms: the lattice
of congruences on a coop is isomorphic to its lattice of ideals in the sense
defined above. We have the following immediate corollary:
\begin{Corollary}\label{cor:coop-simple}
A coop is simple iff its $(0, +, \rImp)$-reduct is a simple hoop. \Done
\end{Corollary}

In categorical language, the forgetful functor from the category of coops to
the category of hoops provides an isomorphism between the category of coops
and the full subcategory of the category of hoops comprising
the objects satisfying the axiom $\all{x}\ex{y}x = y \rImp x$.
In fact, there is a functor that
maps a hoop to an enveloping coop. This is adjoint to the forgetful functor
from coops to hoops. The forgetful functor is faitfhul (as they always are) and
the above says that it is full as well.

\Subsection{Simple Coops}

The hoop $\VH$ is said to be {\em archimedean} iff, for any non-zero
$x \in H$ and any $y \in H$, there is $m \in \NN$, such
that $y \le mx$. We then have:
\begin{Theorem}\label{thm:simple-hoops-archimedean}
A hoop is simple iff it is archimedean.
\end{Theorem}
\Proof
Immediate from the definition of $\I(x)$
and the fact that $\VH$ is simple iff $\I(x) = H$ for every
non-zero $x \in H$.
\Done

\begin{Theorem}\label{thm:simple-coops-archimedean}
A coop is simple iff it is archimedean.
\end{Theorem}
\Proof
Immediate from Corollary~\ref{cor:coop-simple}
and Theorem~\ref{thm:simple-hoops-archimedean}.
\Done

\vspace{4mm}

We will need an interesting property of hoops due to
Bosbach~\cite{bosbach69a}. From a logical perspective, this says that
$\LLi$ enjoys the principle that to prove an implication one may assume the
converse implication.
\begin{Lemma}\label{lma:bosbach}
Let $\VH$ be a hoop, $x, y \in H$. Then
\[
(x \rImp y) \rImp (y \rImp x) = y \rImp x
\]
\end{Lemma}
\Proof
Clearly $y \rImp x \ge (x \rImp y) \rImp y \rImp x$, so
it is enough to prove that $((x \rImp y) \rImp y \rImp x \ge y \rImp x$,
or equivalently that $y + ((x \rImp y) \rImp y \rImp x) \ge x$,
but we have:
\begin{align*}
y + ((x \rImp y) \rImp y \rImp x) &= \\
   y + (y \rImp (x \rImp y) \rImp x) &=
    & \tag*{$\Cwc$}\\
  ((x \rImp y) \rImp x) + (((x \rImp y) \rImp x) \rImp y) &\ge \\
   ((x \rImp y) \rImp x) + (x \rImp y) &\ge x
\end{align*}
\noindent where the penultimate inequality holds since $\rImp$ is antimonotonic in its first argument and $(x \rImp y) \rImp x \le x$.
\Done

\begin{Lemma}\label{lma:m14-linear}
Let $\VH$ be a hoop such that for all $x, y \in H$, if $y = x \rImp y$,
then $x = 0$ or $y = 0$. Then $\VH$ is linearly ordered.
\end{Lemma}
\Proof
By Lemma~\ref{lma:bosbach}, $(a \rImp b) \rImp (b \rImp a) = b \rImp a$
and then by assumption, either $a \rImp b = 0$ or $b \rImp a = 0$, i.e., either
$a \ge b$ or $b \ge a$.
\Done

\begin{Lemma}\label{lma:simple-m14}
If $\VH$ is a simple hoop and $x, y \in H$ are such that $y = x \rImp y$,
then $x = 0$ or $y = 0$.
\end{Lemma}
\Proof
If $y = x \rImp y$, it is easy to see by induction that $y = nx \rImp y$, for
every $n \in \NN$.  But by Theorem~\ref{thm:simple-hoops-archimedean}, $\VH$ is
archimedean, so either $x = 0$ or, for some $n$, $y = nx \rImp y = 0$.
\Done

\begin{Theorem}\label{thm:simple-hoops-linear}
\label{thm:simple-coops-linear}
Simple hoops and simple coops are linearly ordered.
\end{Theorem}
\Proof
For hoops, this is immediate from Lemmas~\ref{lma:simple-m14} and~\ref{lma:m14-linear}.
The statement for coops follows using Corollary~\ref{cor:coop-simple}.
\Done

\vspace{4mm}

We will see later that  simple coops are also Wajsberg hoops.

\begin{Lemma}\label{lma:subcoops-of-real-coops}
Let $\VC$ be a coop such that $C \subseteq \RR$ and let $\VG$ be the subgroup
of the additive group $\RR$ generated by $C$. Then $\VG$ is 2-divisible and:\\
{\em(i)} if $\VC$ is a subcoop of $\RR_{{\ge}0}$, then
$G = C \cup -C$ and $\VC = \VG_{{\ge}0}$;\\
{\em(ii)} if $\VC$ is a subcoop of $\RR_{[0, 1]}$ and $1 \in C$, then
$G = \bigcup_{n\in\ZZ} (n + C)$ and $\VC = \VG_{[0, 1]}$.
\end{Lemma}
\Proof
If $g \in G$, $g$ can be written as $i_1x_1 + \ldots i_mx_m$ where $x_j \in C$
and $i_j \in \ZZ$.  But then $g/2 = i_1y_1 + \ldots i_my_m$, where $y_j = x_j/2
\in C$. So $G$ is indeed 2-divisible.
\\
{\em(i):}
It is enough to prove that $G = C \cup -C$, for then
$C = G \cap \RR_{{\ge}0}$ and so $\VC = \VG_{{\ge}0}$.
Since clearly $C \cup -C \subseteq G$, we have only to show
$C \cup -C$ is closed under negation and addition.
Closure under negation is clear. To show closure under addition,
we have to show that if $x, y \in C$,
then $x + y$, $-x + -y$ and $x - y$ are in $C \cup -C$. This is clear
for $x + y$ and $-x + -y$, since $C$ is closed under addition. As for
$x - y$, if $x \ge y$, then, by definition, $y \rImp x = x - y \in C$,
while, if $x < y$, $x \rImp y = y - x \in C$ and so $x - y \in -C$.\\
{\em(ii):}
It is enough to prove that $G = \bigcup_{n\in\ZZ} (n + C)$, for then
$C = G \cap [0, 1]$ and so $\VC = \VG_{[0, 1]}$.
Clearly $\bigcup_{n\in\ZZ} (n + C) \subseteq G$, so we have only to show that
$\bigcup{n\in\ZZ} (n + c)$ is closed under negation and addition. So let $x, y
\in C$ and $j, k \in \ZZ$ be given. We have:
\[
-(j + x) = -(j+1) + 1 - x = -(j+1) + (x \rImp 1) \in -(j+1) + C
\]
\noindent giving closure under negation. If $x + y \le 1$ (in $\VG$, not
$\VC$), then we have:
\[
(j + x) + (k + y) = (j+k) + (x + y) \in (j+k) + C,
\]
\noindent while if $1 < x + y < 2$, we can find $i, n \in \NN$ with
$i \le 2^n$, such that $x > \frac{i}{2^n}$ and $y > \frac{2^n-i}{2^n}$
and then we have:
\begin{align*}
(j + x) + (k + y) &=
   (j+k+1) + (x - \frac{i}{2^n}) + (y - \frac{2^n-i}{2^n}) \\
   &= (j+k+1) + (\frac{i}{2^n} \rImp x) + (\frac{2^n-i}{2^n} \rImp y) \\
   &\in (j+k+1) + C
\end{align*}
\noindent since $1 \in C$, so that $\frac{i}{2^n}, \frac{2^n-i}{2^n} \in C$,
since $C$ is closed under halving and coop addition (which agrees with the
group addition when the sum in the group is at most 1).
Finally if  $x + y = 2$, we have:
\[
(j + x) + (k + y) = (j+k+2) + 0 \in (j+k+2) + C.
\]
In all cases, $(j + x) + (k + y) \in \bigcup_{n\in\ZZ} (n + C)$
and so $\bigcup_{n\in\ZZ} (n + C)$ is closed under addition, as claimed.
\Done

\vspace{4mm}

Dyadic rational numbers will play an important r\^{o}le in the sequel
as they did in the above proof. We will now generalise the halving operator
on a coop to multiplication by arbitrary non-negative dyadic rationals.
So, let $\VC = (C; 0, +, \rImp, /2)$ be any coop and
define a function $\phi : \NN_{{>}0} \times \NN \times C \To C$ such that:
\begin{align*}
\phi(1, 0, x) &= x \\
\phi(1, n+1, x) &= \phi(1, n, x)/2 \\
\phi(i, n, x) &= i\phi(1, n, x)
\end{align*}
Using the fact that $x/2 + x/2 = x$, we find that the following
holds for any $i, n \in \NN$ and $x \in C$.
\begin{align*}
\phi(2i, n+1, x) &= \phi(i, n, x) \\
\end{align*}
Thus, if $\frac{i}{2^n} = \frac{j}{2^m}$ (in $\QQ$),
$\phi(i, n, x) = \phi(j, m, x)$ for any $x$, and so
$\phi$ induces a function $\Dnn \times C \To C$
which we write multiplicatively: $(p, x) \mapsto px$.
(Here, as with $\NN, \ZZ$, etc., we abuse notation by writing $\DD$, $\Dnn$ and
$\DD_{[0, a]}$ both for the structures and for their carrier sets.)
So for example $\frac{3}{4}x = (x/2)/2 + (x/2)/2 + (x/2)/2$.

Clearly we have $(p+q)x = px + qx$, so, for fixed $x$,
$p \mapsto px$ defines a homomorphism of monoids from $\Dnn$ to $\VC$.
Also, we have $p(x + y) = px + py$, so that for fixed $p$,
$x \mapsto px$ is a homomorphism of monoids from $\VC$ to itself.
If $p, q \in \DD$ with $0 \le p, q \le 1$, we have $p(qx) = (pq)x$, so we have
an action on $\VC$ {\it qua} monoid
of the multiplicative monoid of dyadic rationals in
the interval $[0, 1]$.
However, if $p > 1$ or $q > 1$, $p(qx) \not= (pq)x$ in
general; e.g. with $M = \DD_{[0, 1]}$ and $x = 1$, one has $2x = x$, so that
$\frac{1}{2}(2x) = \frac{1}{2}x = \frac{1}{2}$, while
$(\frac{1}{2}.2)x = 1x = 1$.

\begin{Lemma}\label{lma:dyadic-fractions-in-a-coop}
Let $x \not= 0$ be an element of a coop, $\VC$,
and $0 \le i < j \le 2^n$. Then
{\em(i)} $\frac{i}{2^n}x < \frac{j}{2^n}x$ and
{\em(ii)} $\frac{i}{2^n}x \rImp \frac{j}{2^n}x = \frac{j-i}{2^n}x$.
\end{Lemma}
\Proof We prove {\em(ii)} first. Note that
since $\frac{i}{2^n}x + \frac{j-i}{2^n}x =  \frac{j}{2^n}x$,
we have
$\frac{j-i}{2^n}x \ge \frac{i}{2^n}x \rImp \frac{j}{2^n}x$
by the residuation property.
Thus as $a \rImp b \ge a + c \rImp b + c$,
it is enough to prove {\em(ii)} in the special
case when $j = 2^n$, for then for $j < 2^n$ we have:
\[
\frac{i}{2^n}x \rImp \frac{j}{2^n}x
\ge \frac{i + 2^n - j}{2^n}x \rImp \frac{2^n}{2^n}x
= \frac{2^n - (i + 2^n - j)}{2^n}x
= \frac{j-i}{2^n}x.
\]

So taking $j = 2^n$, let us prove {\em(ii)} by induction on $n$.
The statement is trivial when $n=0$.
So given $n \ge 0$ assume that
$\frac{i}{2^n}x \rImp x = \frac{2^n-i}{2^n}x$ holds for any $x$
and $i$ with $0 \le i < 2^n$. Let $x$ and $i$ with $0 \le i < 2^{n+1}$
be given. If $i = 2^n$, then
$\frac{i}{2^{n+1}}x = \frac{1}{2}x$ and we
have $\frac{1}{2}x \rImp x = \frac{1}{2}x$ by the coop laws.
If $i < 2^n$, we have (using the inductive hypothesis on the line marked ($*$)):
\begin{align*}
\frac{2^{n+1}-i}{2^{n+1}} x
   &= \frac{2^n-i}{2^{n+1}}x + \frac{1}{2}x \\
   &= \frac{2^n-i}{2^{n+1}}x + \left(\frac{1}{2}x \rImp x\right) &\tag*{\HALVE}\\
   &= \frac{2^n-i}{2^{n+1}}x + \left(\frac{2^n-i}{2^{n+1}}x  + \frac{i}{2^{n+1}}x \rImp x\right) \\
   &= \frac{2^n-i}{2^{n+1}}x + \left(\frac{2^n-i}{2^{n+1}}x  \rImp \frac{i}{2^{n+1}}x \rImp x\right) \\
   &= \left(\frac{i}{2^{n+1}}x \rImp x\right) + \left[\left(\frac{i}{2^{n+1}}x \rImp x\right) \rImp \frac{2^n-i}{2^{n+1}}x\right] &\tag*{$\Cwc$}\\
   &= \left(\frac{i}{2^{n+1}}x \rImp x\right) + \left[\left(\frac{i}{2^{n+1}}x \rImp x\right) \rImp \frac{2^n-i}{2^{n}}\frac{1}{2}x\right] \\
   &= \left(\frac{i}{2^{n+1}}x \rImp x\right) + \left[\left(\frac{i}{2^{n+1}}x \rImp x\right) \rImp \left(\frac{i}{2^n}\frac{1}{2}x \rImp \frac{1}{2}x\right)\right] &\tag*{($*$)} \\
   &= \left(\frac{i}{2^{n+1}}x \rImp x\right) + \left[\left(\frac{i}{2^{n+1}}x \rImp x\right) \rImp \left(\frac{i}{2^{n+1}}x \rImp \frac{1}{2}x\right)\right] \\
   &= \frac{i}{2^{n+1}}x \rImp x.
\end{align*}
If $2^{n+1} > i > 2^n$, then we have:
\begin{align*}
\frac{i}{2^{n+1}}x \rImp x 
   &= \frac{i-2^n}{2^{n+1}}x \rImp \frac{1}{2}x \rImp x \\
   &= \frac{i-2^n}{2^{n+1}}x \rImp \frac{1}{2}x &\tag*{\HALVE}\\
   &= \frac{i-2^n}{2^n}\frac{1}{2}x \rImp \frac{1}{2}x \\
   &= \frac{2^n-(i-2^n)}{2^n}\frac{1}{2}x &\tag*{($*$)}\\
   &= \frac{2^{n+1}-i}{2^{n+1}}x.
\end{align*}
\noindent
This completes the proof of part {\em(ii)}.
Part {\em(i)} follows since, by part {\em(ii)},
we have $\frac{i}{2^n}x \rImp \frac{i+1}{2^n}x = \frac{1}{2^n}x \not= 0$,
whence $\frac{i}{2^n}x < \frac{i+1}{2^n}x \le \frac{j}{2^n}x$.
\Done

\vspace{4mm}

By the following lemma, simple coops are semi-cancellative.

\begin{Lemma}\label{lma:simple-semi-cancellative}
Let $\VC$ be a coop and let $x, y \in C$ be such that $x + y = x$,
then either $x = mx$ for all $m \in \NN$ or $y \le \frac{1}{2^n}x$ for
all $n \in \NN$.
In particular, if $\VC$ is simple, and hence archimedean, either
$x$ is an annihilator or $y = 0$.
\end{Lemma}
\Proof
By an easy induction, we have $x + my = x$ for all $m \in \NN$.
If $y > \frac{1}{2^n}x$ for some $n \in \NN$, then we have
$$
x = x + 2^ny \ge x + 2^n(\frac{1}{2^n}x) = 2x
$$
\noindent
and then by another easy induction we have $x = mx$ for all $m \in \NN$.
\Done

\vspace{4mm}

If $\VH$ is a hoop and $0 \not= x \in H$, define the {\em depth} of $x$ to be
the smallest $d \in \NN$ such that $(d+1)x = dx$, or to be $\infty$ if no such
$d$ exists.  Lemma~\ref{lma:dyadic-fractions-in-a-coop} implies that if $x$ is
a non-zero element of a coop, then the depth of $\frac{1}{2^n}x$ is at least
$2^n$.

\begin{Lemma}\label{lma:simple-hoop-depth}
Let $\VH$ be a simple hoop. Then either
{\em(i)} every non-zero element has infinite depth or
{\em(ii)} $\VH$ is bounded and every non-zero element has finite depth.
\end{Lemma}
\Proof
Assume {\em(i)} does not hold, so there is a non-zero
$x \in H$ with finite depth $d$, so $dx = (d+1)x$.
By induction, for any $n > d$, we have $dx = nx$.
Let $a = dx$. Then $na = ndx = dx = a$ for any $n > 0$, so
that, as $\VH$ is simple,  $H = \I(a) = \Downset{a}$.
Now if $y$ is any non-zero element, $\I(y) = H$, so $a \le ny$ for some $n$
and we have $ny \ge a \ge (n+1)y$ so that $ny = (n+1)y$ and $y$ has finite
depth.
\Done

\begin{Theorem}\label{thm:simple-coop-structure}
Let $\VC$ be a simple coop.
Then there is 2-divisible subgroup $\VG$ of the additive group $\RR^{+}$,
such that either
{\em(i)} $\VC$ is isomorphic to $\VG_{{\ge}0}$, or
{\em(ii)} $\VC$ is isomorphic to $\VG_{[0, 1]}$.
\end{Theorem}
\Proof
By Theorems~\ref{thm:simple-coops-archimedean} and~\ref{thm:simple-coops-linear},
$\VC$ is archimedean and linearly ordered. We use these properties without further
comment in the rest of the proof. \\
If $C$ is not bounded, then, by Lemma~\ref{lma:simple-hoop-depth}, there is
a non-zero $e \in C$ with infinite depth, so that
$ne < (n+1)e$ for every $n \in \NN$. We will show that case {\em(i)} holds.
To see this, define $f : C \To \RR_{{\ge}0}$ by:
\[
f(x) = \Sup \{\frac{i}{2^n} \ST i, n \in \NN, \frac{i}{2^n}e \le x\}.
\]
For every $x$, we have $0e \le x \le ne$ for all large enough $n$.
Thus the set whose supremum is used in the definition of $f$ is
always non-empty and bounded above, so $f$ is well-defined.
Clearly, $f$ is at least weakly monotonic, i.e., if $x \le y$, then $f(x) \le f(y)$.
We claim that $f$ is a homomorphism from $\VC$ to the coop $\RR_{{\ge}0}$.
By Lemma~\ref{lma:simple-semi-cancellative}, the $(0, +)$-reduct of
$\VC$ is a cancellation
monoid, so that as $2e = e + e = e + (e \rImp 2e)$, we have $e = e \rImp 2e$,
whence $e = (2e)/2$. By induction, $e = \frac{1}{2^n}(2^ne)$ for any
$n \in \NN$. Hence, for any $m$, taking $x = 2^{m}e$ in
Theorem~\ref{lma:dyadic-fractions-in-a-coop} we find that for any $a, b$ of the
form $\frac{i}{2^n}x$, with $0 \le i \le 2^n$ we have $f(a + b) = f(a) + f(b)$
and $f(a \rImp b) = f(a) \rImp f(b)$. Letting $m$
tend to infinity, these equations hold for any $a, b \in D$, where
$D$ is the set $\{\frac{i}{2^n}e \ST i, n \in \NN\}$ of all dyadic rational multiples
of $e$. Thus $D$ is a subcoop of $\VC$ isomorphic to  $\DD_{{\ge}0}$.
By Theorem~\ref{lma:dyadic-fractions-in-a-coop}, given $x, y \in C$ and any $n \in \NN$,
there are $p_n, q_n, r_n, s_n \in D$ such that $p_n \le x \le q_n$, $r_n \le y \le s_n$,
$f(q_n) - f(p_n) = \frac{1}{2^{n+1}}$ and $f(s_n) - f(r_n) = \frac{1}{2^{n+1}}$. 
But then as $f|_D$ is a coop-homorphism and $f$ is weakly monotonic, we have:
\[
\begin{array}{c}
f(p_n + r_n) = f(p_n) + f(r_n) \\
f(q_n + s_n) = f(q_n) + f(s_n) \\
f(p_n + r_n) \le f(x + y) \le f(q_n + s_n) \\
f(p_n) + f(r_n) \le f(x) + f(y) \le f(q_n) + f(s_n) \\
f(q_n + s_n) - f(p_n + r_n) \le \frac{1}{2^n} \\
\\
f(q_n \rImp r_n) = f(q_n) \rImp f(r_n) \\
f(p_n \rImp s_n) = f(p_n) \rImp f(s_n)  \\
f(q_n \rImp r_n) \le f(x \rImp y) \le f(p_n \rImp s_n) \\
f(q_n) \rImp f(r_n) \le f(x) \rImp f(y) \le f(p_n) \rImp f(s_n) \\
f(p_n \rImp s_n) - f(q_n \rImp s_n) \le \frac{1}{2^n}
\end{array}
\]
\noindent
Letting $n$ tend to infinity, we must have that $f(x + y) = f(x) + f(y)$ and $f(x \rImp y) = f(x) \rImp f(y)$
and $f$ is indeed a homomorphism from $\VC$ to $\RR_{{\ge}0}$ as claimed.
But $\VC$ is simple, hence $f$ is either identically zero or is one-to-one, but
clearly $f(e) = 1 \not= 0$, so $f$ embeds $\VC$ as a subcoop of $\RR_{{\ge}0}$.
{\em(i)} follows immediately using Lemma~\ref{lma:subcoops-of-real-coops}. \\
Now assume $\VC$ is bounded, with annihilator $a$, say. So $a \ge x$ for every
$x \in C$.  To see that that case {\em(ii)} holds, define $g : C \To
\RR_{[0,1]}$ by:
\[
g(x) = \Sup \{\frac{i}{2^n} \ST i, n \in \NN, i \le 2^n, \frac{i}{2^n}a \le x\}.
\]
Then by an argument very similar to the one used above in  the unbounded case,
$C$ has a dense subcoop $D_1$ such that $g|_{D_1}$ is an
isomorphism from $D_1$ to $\DD_{[0,1]}$. Then, approximating $x + y$ and $x
\rImp y$ by elements of $D_1$ just as we did above, we find that $g$ is a
homomorphism embedding $\VC$ as a subcoop of $\RR_{[0,1]}$, from which
{\em(ii)} follows using Lemma~\ref{lma:subcoops-of-real-coops}.
\Done

\Subsection{Subdirectly Irreducible Coops}
Recall that an algebra $\VA$ is {\em subdirectly irreducible} iff
the intersection $\mu$ of its non-identity congruences is not the identity
congruence, in which case $\mu$ is called the {\em monolith}.
Thus a hoop or a coop is subdirectly irreducible iff the
intersection of all its non-zero ideals is non-zero. In this section,
we determine the structure of subdirectly irreducible coops.

If $\VC$ and $\VD$ are subcoops of a coop $\VE$, we say
$\VE$ is the {\em ordinal sum} of $\VC$ and $\VD$
and write $\VE = \VC \ordSum \VD$ iff
$C \cap D = \{0\}$, $C \cup D = E$ and whenever
$c \in C$ and $0 \not= d \in D$, $c + d = d$.
It is easy to see that, if $\VE = \VC \ordSum \VD$ and $c \in C$ and $0 \not= d \in D$,
then $d > c$ and $c \rImp d = d$.
Thus $C$ is an ideal and $\VE/C \cong \VD$.
We will find that any subdirectly irreducible coop is $\VS \ordSum \VF$
where $\VS$ is totally ordered and subdirectly irreducible and $\VF$ can be any coop.
This could also be establishing
using the analogous result for subdirectly irreducible hoops
proved in \cite{blok-ferreirim00}, but the extra structure in a coop admits
a slightly more efficient presentation.

\begin{Theorem}\label{thm:cep}
Hoops and coops have the congruence extension property.
\end{Theorem}
\Proof
By Theorem~\ref{thm:coop-homomorphisms} and the discussion of ideals
that precedes it, it suffices to show that if $\VC$ is a subhoop of
a hoop $\VD$, then for any ideal $I \subseteq C$, there is an ideal
$J \subseteq D$, such that $I = J \cap C$. But, if $I \subseteq C$
is an ideal, it is easily verified from the definitions that $I = J \cap C$
where $J$ is the ideal of $\VD$ generated by $I$.
\Done

\vspace{4mm}

Let $\VC$ be a subdirectly irreducible coop, so that the set of all
non-zero ideals of $\VC$ intersect in a non-zero ideal $M$, which
we call the {\em monolithic ideal}. Since coops have the congruence extension
property, $M$ viewed as a coop in its own right can have no non-trivial ideals,
so $M$ is a simple coop, and so by Theorems~\ref{thm:simple-coops-archimedean}
and~\ref{thm:simple-coops-linear}, $M$ is archimedean and linearly ordered.

If $x \in C$, we define the {\em implicative stabilizer} $\IS(x)$ as follows:
\begin{align*}
\IS(x) &\IsDef \{s \in C \ST  s \rImp x = x\}.
\end{align*}
\noindent
It is easily verified that $\IS(x)$ is an ideal.
So, for any $x$, either $\IS(x) = \{0\}$ or $\IS(x) \supseteq M$.
If $X \subseteq C$, we write $\IS(X)$ for $\bigcap_{x\in X} \IS(x)$.

\begin{Theorem}\label{thm:subdirectly-irreducible-coops}
Let $\VC$ be a subdirectly irreducible coop with monolithic ideal $M$ and
let $F, S \subseteq C$ be defined as follows:
\begin{align*}
F &\IsDef \{f \in C \ST M \subseteq \IS(f)\} \\
S &\IsDef \IS(F)
\end{align*}
Then:
\begin{align*}
\mbox{{\em(i)}} \quad &
   \all{x \in C \Diff \{0\}}\ex{a \in M \Diff \{0\}} x \ge a; \\
\mbox{{\em(ii)}} \quad &
   \all{f \in F \Diff \{0\}, a \in M} f \ge a; \\
\mbox{{\em(iii)}} \quad &
   \all{a \in M, f \in F \Diff \{0\}} a + f = f; \\
\mbox{{\em(iv)}} \quad &
   \all{f \in F \Diff \{0\}, x \in C} x \ge f \Imp x \in F; \\
\mbox{{\em(v)}} \quad &
   \all{x \in C, f \in F} x \rImp f \in F; \\
\mbox{{\em(vi)}} \quad &
   \all{f \in F \Diff \{0\}, x \in C \Diff F} f > x; \\
\mbox{{\em(vii)}} \quad &
   \mbox{$F$ is the carrier set of a subcoop $\VF$ of $\VC$}; \\
\mbox{{\em(viii)}} \quad &
   \mbox{$S$ is a linearly ordered ideal of $\VC$, and $S \cap F = \{0\}$}; \\
\mbox{{\em(ix)}} \quad &
  \mbox{Writing $\VS$ for the subcoop with carrier set $S$,
          $\VS$ is semi-cancellative}; \\
\mbox{{\em(x)}} \quad &
   \VC = \VS \ordSum \VF.
\end{align*}
\end{Theorem}
\Proof
First note that if $x \in C$ and $a \rImp x = x$ for some $a \in C \Diff \{0\}$, then
$\IS(x) \not= \{0\}$, hence $M \subseteq \IS(x)$ so that $x \in F$. \\
{\em(i):}
if $0 \not= x \in C$, then as $\{0\} \not= M \subseteq \I(x)$, there is $a \in M$ and
$n \in \NN$, with $2^{n}x \ge a \not= 0$, but then $0 \not= \frac{1}{2^n}a \in M$
and $x \ge \frac{1}{2^n}(2^{n}x) \ge \frac{1}{2^n}a$ (where the first
inequality follows by induction using  part {\em(iv)} of
Corollary~\ref{cor:halving-miscellany}). \\
{\em(ii):}
since $a \in M \subseteq \I(f)$, $mf \ge a$ for some $m \in \NN$.
By Lemma~\ref{lma:bosbach}, $f \rImp a = (a \rImp f) \rImp f \rImp a = f \rImp f \rImp a = 2f \rImp a$,
since $f \in F$. By induction, $f \rImp a = nf \rImp a$ for every $n \in \NN$.
In particular, $f \rImp a = mf \rImp a = 0$. \\
{\em(iii):}
by part {\em(ii)}, $f \ge a$, i.e. $f \rImp a = 0$. Hence, using $\Cwc$,
$a + f = a + (a \rImp f) = f + (f \rImp a) = f$. \\
{\em(iv)}
assume $f \in F$, $x \in C$ and $x \ge f \not= 0$.
We need to show that if $a \in M$,
$a \rImp x = x$. But given $a \in M$, we have $a \rImp x \ge a \rImp f = f$,
i.e., $((a \rImp x) \rImp f) = 0$.
Hence:
\begin{align*}
a \rImp x 
   &= (a \rImp x) + ((a \rImp x) \rImp f) \\
   &= f + (f \rImp a \rImp x) \tag*{$\Cwc$}\\
   &= f + (a \rImp f \rImp x) \\
   &= f + a + (a \rImp f \rImp x) \tag*{(iii)} \\
   &= f + (f \rImp x) + ((f \rImp x) \rImp a) \tag*{$\Cwc$} \\
   &\ge x.
\end{align*}
So $x \ge a \rImp x \ge x$ giving $x = a \rImp x$ as required.\\
{\em(v):}
if $a \in M$ and $f \in F$, $a \rImp f = f$ by the definition of $F$.
So, for any $x \in C$, $a \rImp x \rImp f = x \rImp a \rImp f = x \rImp f$,
whence $x \rImp f \in F$. \\
{\em(vi):}
Let $f \in F$ and $x \in C \Diff F$.
if $a \in M$, we have $a \rImp f \rImp x = a + f \rImp x = f \rImp x$, by part
{\em(iii)}, so $x \ge f \rImp x \in F$ and by part {\em(iv)} we can only have
$f \rImp x = 0$, i.e., $f \ge x$, and the inequality must be strict,
since  $x \not\in F$. \\
{\em(vii):} Clearly $0 \in F$. Given $f, g \in F$,
we must show that $f + g, f \rImp g$ and $f/2$ all belong to $F$.
As $f + g \ge f$, $f + g \in F$ follows from part {\em(iv)}.
That $f \rImp g \in F$ follows from part {\em(v)}.
Finally $f/2 \in F$ follows from part {\em(iii)} together
with part {\em(v)} of Corollary~\ref{cor:halving-miscellany},
since given $0 \not= a \in M$ and $a \rImp f = f$, then we
have $0 \not= a/2 \in M$ and $a/2 \rImp f/2 = (a \rImp f)/2 = f/2$.\\
{\em(viii):} As the intersection of a set of ideals, $S$ is itself an ideal.
If $x \in S \cap F$, then $x = x \rImp x = 0$, so $S \cap F = \{0\}$.
If $s, t \in S$ and $s \rImp t = t$, I claim that either $s = 0$
or $t = 0$, whence $S$ is linearly ordered by Lemma~\ref{lma:m14-linear}.
To prove the claim, if $s \rImp t = t$ and $s \not= 0$, then, by part {\em(i)},
there is $a \in M$ such that $s \ge a > 0$, but then
$t \ge a \rImp t \ge s \rImp t \ge t$, so $a \rImp t = t$ and $t \in F$,
so $t \in S \cap F$ and therefore $t = 0$. \\
{\em(ix):}
let $s, t \in S$ with $t \not= 0$ and $s + t = s$.
We must prove that $s$ annihilates $S$, i.e., $S = \Downset{s}$.
We may choose an $a \in M$ such that $t \ge a \not= 0$ and then
$s = s + t \ge s + a \ge s$, whence $s + a = s$.
If $u \in S$, then we have $a \rImp s \rImp u = a + s \rImp u = s \rImp u$,
so $s \rImp u \in S \cap F = \{0\}$ by part {\em(viii)}.
Hence $s \rImp u = 0$, i.e., $s \ge u$. \\
{\em(x):}
by part {\em(viii)} $S \cap F = \{0\}$.
I claim that, if $x \in C \Diff F$ and $0 \not= f \in F$,
then $x \in S$ and $x + f = f$.
Given this, we must have that $C = S \cup F$ and $S \cap F = \{0\}$
and so $\VC = \VS \ordSum \VF$.
So assume $x \in C \Diff F$ and $0 \not= f \in F$.
We must prove that $f = x \rImp f = x + f$.
By part {\em(v)}, $(x \rImp f) \rImp f \in F$, but
$x \ge (x \rImp f) \rImp f$ and $x \not\in F$, so by part {\em(iv)},
we must have $(x \rImp f) \rImp f = 0$, i.e., $x \rImp f \ge f$ implying
$f = x \rImp f$. Using $\Cwc$ and part {\em(vi)}, we have
$f = f + (f \rImp x) = x + (x \rImp f) = x + f$ and the claim is true.
\Done

\vspace{4mm}

We refer to the subcoops $\VF$ and $\VS$ of the theorem
as the {\em fixed} subcoop and the {\em support} subcoop respectively.
Since the support is linearly ordered and semi-cancellative
the following theorem applies to it.

\begin{Theorem}\label{thm:lin-semi-canc}
Let $\VL$ be a linearly ordered semi-cancellative coop. Then
\\
{\em(i)}
If $\VL$ is bounded, it is involutive; \\
{\em(ii)}
$\VL$ is a Wajsberg coop, i.e.,
for any $s, t \in L$, 
$ (t \rImp s) \rImp s = (s \rImp t) \rImp t $.
\end{Theorem}
\Proof
{\em(i):}
As usual write $1$ for the annihilator of $\VL$ and $\Not x$ for $x \rImp 1$.
Note that by Corollary \ref{cor:halving-miscellany},
$(\Not s)/2 = s/2 \rImp 1/2$.
We must show that $\Not\Not s = s$ for any $s \in L$.
We claim that $\Not(\cdot) : L \To L$ is injective.
To see this, assume $s, t \in L$ with $\Not s = \Not t$.
Then $s/2 \rImp 1/2 = (\Not s)/2 = (\Not t)/2 = t/2 \rImp 1/2$.
As $1/2 \ge s/2$ and $1/2 \ge t/2$,
we have $1/2 = (\Not s)/2 + s/2 = (\Not t)/2 + t/2$.
But $1/2$ is not an annihilator, so by part {\em(i)}, this implies
$s/2 = t/2$, whence $s = t$, completing the proof that $\Not$ 
is injective.
But for any $s \in L$, $\Not\Not\Not s = \Not s$,
so if $\Not$ is injective, $\Not\Not s = s$. \\
{\em(ii):}
\noindent
We claim that for any $s, t \in L$, $(t \rImp s) \rImp s = \Min\{s, t\}$,
which is well-defined because $\VL$ is linearly ordered.
Assuming the claim, we have:
$$
(t \rImp s) \rImp s = \Min\{s, t\} = \Min\{t, s\} = (s \rImp t) \rImp t
$$
\noindent
so the claim implies the required identity.
To prove the claim, note that if $t \ge s$, we have:
$$
(t \rImp s) \rImp s = 0 \rImp s = s = \Min\{s, t\}
$$
\noindent
while if $s \ge t$, we have
$$
s = (t \rImp s) + t = (t \rImp s) + ((t \rImp s) \rImp s).
$$
If $\VL$ is unbounded or if $\VL$ is bounded but $s$ is not
the annihilator, then the semi-cancellative property gives us
$$
(t \rImp s) \rImp s = t = \Min\{s, t\}.
$$
\noindent
If $\VL$ is bounded, let us write $1$ for its annihilator and $\Not x$ for
$x \rImp 1$ as we did in the proof of part {\em(i)}. Then
if $s = 1$, we have:
$$
(t \rImp s) \rImp s = \Not\Not t = t = \Min\{s, t\}
$$
\noindent
by part {\em(i)}.
In all cases, the claim holds and the proof is complete.
\Done

\vspace{4mm}

Note that a non-trivial ordinal sum is never Wajsberg: if $s \in S$
and $f \in F$, then in $\VS \ordSum \VF$, we have:
\begin{align*}
(s \rImp f) \rImp f &= f \rImp f &= 0 \\
(f \rImp s) \rImp s &= 0 \rImp s &= s
\end{align*}
So $(s \rImp f) \rImp f = (f \rImp s) \rImp s$ iff $0 \in \{s, f\}$.

\begin{Theorem}\label{thm:wajsberg-coops-univ-decidable}
The universal theory of Wajsberg coops is decidable.
(I.e., the set of purely universal formulas in the language of
a coop that are valid in all coops is decidable).
\end{Theorem}
\Proof
We claim that any Wajsberg coop is isomorphic to a subcoop of a product of
linearly ordered semi-cancellative coops.
Given part {\em(ii)} of Theorem~\ref{thm:lin-semi-canc}, such a product
is itself a Wajsberg hoop, hence, given the claim, the universal theory of Wajsberg coops
reduces to that of linearly ordered semi-cancellative coops and by
Theorem~\ref{thm:lin-canc-coops-decidable} the full first order theory
of linearly ordered semi-cancellative coops is decidable. \\
As for the claim, let $\VW$ be a Wajsberg coop. By Birkhoff's theorem,
$\VW$ embeds in a product $\prod_i \VC_i$, where each $\VC_i$
is a subdirectly irreducible homomorphic image of $\VW$.
By the remarks above, when we write $\VC_i$ as the ordinal sum
of its support and fixed part, $\VS_i \ordSum \VF_i$, 
$\VF_i = \{0\}$, so $\VC_i$ is isomorphic to $\VS_i$.
The claim follows from parts {\em(viii)} and {\em(ix)}
of Theorem~\ref{thm:subdirectly-irreducible-coops} and
Theorem~\ref{thm:lin-semi-canc}.
\Done

\Section{Future Work}\label{sec:future-work}
An important goal of our work is to understand the decision problem
for useful classes of coop, and we have presented some results
on Wajsberg coops, in particular, in the present paper.
We already have some more results about general coops,
but the proofs are not yet in a very satisfactory form.
Blok and Ferreirim have shown that the quasi-equational theory
of hoops is decidable. Using their results on subdirectly
irreducible hoops one can show that any hoop embeds in a coop and
from this conclude that the quasi-equational theory of coops
is decidable (since this implies that a horn clause in the
language of coops, can be translated into an equisatisfiable
Horn clause in the language of hoops\footnote{
To do this, replace subterms of the form $t/2$ by $v_t$ where $v_t$
is a fresh variable and add hypotheses $v_t = v_t \rImp t$.}
However, this approach doesn't yield a practically feasible algorithm
and gives no information about the complexity of the decision problem.
We hope to improve on this position in future work.

\bibliographystyle{plain}
\bibliography{coops}

\begin{thebibliography}{10}

\bibitem{ben-yaacov08}
I.~{Ben Yaacov}.
\newblock {On theories of random variables}.
\newblock Available on line at: {\tt http://arxiv.org/TBS.TBS}, 2008.

\bibitem{ben-yaacov-pedersen09}
I.~{Ben Yaacov} and A.~P. {Pedersen}.
\newblock {A proof of completeness for continuous first-order logic}.
\newblock Available on line at: {\tt http://arxiv.org/0903.4051}, 2009.

\bibitem{blok-ferreirim00}
W.~J. Blok and I.~M.~A. Ferreirim.
\newblock On the structure of hoops.
\newblock {\em Algebra Universalis}, 43(2-3):233--257, 2000.

\bibitem{bosbach69a}
B.~Bosbach.
\newblock {Komplement\"are Halbgruppen. Axiomatik und Arithmetik.}
\newblock {\em Fundam. Math.}, 64:257--287, 1969.

\bibitem{chang58b}
C.~C. Chang.
\newblock Algebraic analysis of many valued logics.
\newblock {\em Trans. Amer. Math. Soc.}, 88:467--490, 1958.

\bibitem{chang58a}
C.C. Chang.
\newblock {Proof of an axiom of {\L}ukasiewicz.}
\newblock {\em Trans. Am. Math. Soc.}, 87:55--56, 1958.

\bibitem{Ciabattoni:1997}
Agata Ciabattoni and Duccio Luchi.
\newblock Two connections between linear logic and {L}ukasiewicz logics.
\newblock In {\em Proceedings of the 5th Kurt G\"odel Colloquium on
  Computational Logic and Proof Theory}, pages 128--139, London, UK, 1997.
  Springer-Verlag.

\bibitem{Girard(87B)}
J.-Y. Girard.
\newblock Linear logic.
\newblock {\em Theoretical Computer Science}, 50(1):1--102, 1987.

\bibitem{Hajek98}
Petr H{\'a}jek.
\newblock {\em Metamathematics of Fuzzy Logic}.
\newblock Kluwer Academic Publishers, 1998.

\bibitem{henson-iovino02}
C.~Ward Henson and Jos{\'{e}} Iovino.
\newblock Ultraproducts in analysis.
\newblock In {\em Analysis and Logic}, volume 262 of {\em London Mathematical
  Society Lecture Notes}, pages 1--113. Cambridge University Press, 2002.

\bibitem{Koehler81}
Peter Koehler.
\newblock {Brouwerian semilattices.}
\newblock {\em Trans. Am. Math. Soc.}, 268:103--126, 1981.

\bibitem{lukasiewicz-tarski30}
J.~{\L}ukasiewicz and A.~Tarski.
\newblock {Untersuchungen \"uber den Aussagenkalk\"ul.}
\newblock {\em C. R. Soc. Sc. Varsovie 23}, (1930):30--50, 1930.

\bibitem{raftery07}
James~G. Raftery.
\newblock On the variety generated by involutive pocrims.
\newblock {\em Rep. Math. Logic}, 42:71--86, 2007.

\bibitem{rose-rosser58}
Alan Rose and J.Barkley Rosser.
\newblock {Fragments of many-valued statement calculi.}
\newblock {\em Trans. Am. Math. Soc.}, 87:1--53, 1958.

\end{thebibliography}

\end{document}